\documentstyle[12pt,amssymb,amscd]{amsart}  

\theoremstyle{remark}

\theoremstyle{definition}

\numberwithin{equation}{subsection}

\newcommand{\bbA}{{\Bbb A}}

\newcommand{\bbC}{{\Bbb C}}
\newcommand{\bbD}{{\Bbb D}}

\newcommand{\bbR}{{\Bbb R}}

\newcommand{\bbZ}{{\Bbb Z}}
\newcommand{\cF}{{\cal F}}
\newcommand{\cG}{{\cal G}}
\newcommand{\cO}{{\cal O}}
\newcommand{\cP}{{\cal P}}
\newcommand{\cL}{{\cal L}}
\newcommand{\cM}{{\cal M}}
\newcommand{\cN}{{\cal N}}
\newcommand{\cD}{{\cal D}}
\newcommand{\cV}{{\cal V}}
\newcommand{\cA}{{\cal A}}

\newcommand{\cC}{{\cal C}}
\newcommand{\cE}{{\cal E}}

\newcommand{\cH}{{\cal H}}
\newcommand{\cX}{{\cal X}}
\newcommand{\cY}{{\cal Y}}

\newcommand{\GD}{{\frak D}}
\newcommand{\GG}{{\frak G}}
\newcommand{\Gj}{{\frak j}}

\newcommand{\Ind}{\operatorname{Ind}}  
\newcommand{\Mat}{\operatorname{Mat}}
\newcommand{\Der}{\operatorname{Der}}

\newcommand{\Coker}{\operatorname{Coker}}
\newcommand{\Aut}{\operatorname{Aut}}

\newcommand{\Spf}{\operatorname{Spf}}

\newcommand{\Zar}{\operatorname{Zar}}

\newcommand{\Ker}{\operatorname{Ker}}

\newcommand{\Hom}{\operatorname{Hom}}
\newcommand{\End}{\operatorname{End}}

\newcommand{\isomoto}{\overset{\sim}{\to}}
\newcommand{\id}{\operatorname{Id}}

\newcommand{\shHom}{\underline{\operatorname{Hom}}}

\newcommand{\Conn}{\operatorname{Conn}}

\newcommand{\Smbl}{\operatorname{Smbl}}

\newcommand{\Lie}{\operatorname{Lie}}
\newcommand{\Vect}{\operatorname{Vect}}

\newcommand{\Anch}{\operatorname{Anch}}
\newcommand{\Bun}{\operatorname{Bun}}
\newcommand{\VVect}{{ {\Bbb V}}\hskip -.1cm \operatorname{ect}}
\newcommand{\Spec}{\operatorname{Spec}}
\newcommand{\btimes}{\overline{\otimes}}

\newcommand{\Gg}{{\frak g}}

\newcommand{\gr}{\operatorname{gr}}
\newcommand{\GA}{{\mathfrak A}}
\newcommand{\Ob}{\operatorname{Ob}}
\newcommand{\Mor}{\operatorname{Mor}}
\newcommand{\Res}{\operatorname{Res}}

\newcommand{\Rin}{\operatorname{Rin}}
\newcommand{\FR}{\operatorname{FR}}
\newcommand{\FL}{\operatorname{FL}}

\newcommand{\Mod}{\operatorname{Mod}}
\newcommand{\Rep}{\operatorname{Rep}}
\newcommand{\red}{\operatorname{red}}
\newcommand{\cont}{\operatorname{cont}}

\begin{document}

\title{Free Lie algebroids and the space of paths}

\author{Mikhail Kapranov}
\email{mikhail.kapranov@@yale.edu}
\address{Department of Mathematics, Yale University, 
10 Hillhouse Avenue, New Haven CT 06520 USA}
\maketitle

\vskip 2cm

\centerline {\bf Introduction.}

\vskip 1cm

\noindent {\bf (0.1) } The goal of this paper is to construct algebraic and
algebro-geometric models for  spaces of paths. Paths in a manifold $X$, considered
up to reparametrizations and cancellations, form a groupoid $\Pi_X$.
We construct a Lie algebroid $\cP_X$ which plays the role of the Lie algebra
for $\Pi_X$, i.e., describes ``infinitesimal paths''. 

\vskip .2cm

When $X$ is an algebraic variety, we construct an algebro-geometric model for the
formal neighborhood of $X$ (constant paths) in $\Pi_X$. This is a certain ind-scheme,
denoted $\widehat{\Pi}_X$, which also has a groupoid structure. It is analogous to
(but different from) $\cL X$, the ind-scheme of formal parametrized loops in $X$,
constructed in [KV1]. 

\vskip .2cm

Among known algebro-geometric constructions that have the flavor of the space of unparametrized
paths, one should mention the Kontsevich moduli stack $\overline{M}_{g,2}(X,\beta)$
of stable maps from 2-pointed curves $(C,x,y)$ of genus $g$ into $X$, see [Ma]. 
Indeed, this stack consists of maps $(C,x,y)\to X$ for variable $(C,x,y)$ modulo
isomorphisms of curves, i.e., changes of parametrization. The points $x$ and $y$
serve as the beginning and end of a path. 
One of our main results, Theorem 7.3.5, constructs, in the case $g=0$, a morphism from a certain
formal neighborhood in $\overline{M}_{0,2}(X,\beta)$ into $\widehat{\Pi}_X$. The case of $n$-pointed
curves, $n\geq 3$, can be treated similarly by considering the simplicial classifying space of
$\widehat{\Pi}_X$. 

\vskip .3cm

\noindent {\bf (0.2)} Any vector bundle with connection on $X$ gives a representation of $\Pi_X$
by holonomy. In fact,  representations of $\Pi_X$ are more or less the same as connections, as was shown by Kobayashi [Ko]. 
Therefore it is not so surprising that our Lie algebroid $\cP_X$ can be seen as the
universal receptacle of the curvature data for all connections; these data include
the ``higher covariant derivatives of the curvature''. The quotes here emphasize the
fact that these covariant derivatives are not really defined as tensors; in fact
$\cP_X$ is filtered, with  graded quotients  expressed as
certain tensor spaces. It appears that $\cP_X$ is a fundamental differential-geometric
object  associated to any manifold, similar in importance to the sheaf 
of differential operators. 

\vskip .2cm

Our other main result, Theorem 4.4.3, identifies formal germs of connections on $X$ at $x$
with representations of a certain infinite-dimensional Lie algebra $\cP(X,x)$, which we call
the {\em fundamental Lie algebra} (by analogy with the fundamental group and its relation
with flat connections). Combined with some results of Reutenauer, this provides a
``Taylor formula'' for connections (Corollary 4.4.5). 

\vskip .3cm

\noindent {\bf (0.3)} It is known since the work of K.-T. Chen [C] that the group of paths in 
$\bbR^n$ can be seen as a certain continuous analog of the free group on $n$ generators,
so its Lie algebra is a certain completion of the free Lie algebra.
Our construction of  $\cP_X$
is also a version of the free Lie algebra construction but in the context of Lie algebroids.
This construction generalizes that of Casas, Ladra and Pirashvili [CLP]. 
In the curvature language the appearance of the (free) Lie algebras here
reflects the interpretation of the Bianchi identity  as the
Jacobi identity for the curvature operators $\nabla_i = \partial_i + A_i$. 

\vskip .2cm

In fact, transformation rules for sections of $\cP_X$ from one coordinate system
$(x_1, ..., x_n)$ to another, $(y_1, ..., y_n)$, reflect the transformation rules for 
the system of Chen's iterated integrals of 1-forms constant with respect to these systems:
$$\biggl\{ \int_\gamma^{\rightarrow} dx_{i_1} ... dx_{i_p} \biggr\} \quad \rightsquigarrow \quad
\biggl\{ \int_\gamma^{\rightarrow} dy_{i_1} ... dy_{i_p}\biggr\}.$$
Here the path $\gamma$ is assumed to be the same on both sides. 

\vskip .2cm

Sections of $\cP_X$ give rise to ``noncommutative vector fields'', i.e., to natural systems of
operators $\{P_{E,\nabla}: E\to E\}$ defined for all bundles with connections $(E,\nabla)$ on $X$ and
satisfying the Leibniz rule with respect to the tensor product. By applying the enveloping algebra
construction we get a sheaf of rings $\bbD_X = U(\cP_X)$ whose sections are called
noncommutative differential operators. Such an operator $P$ gives a natural system of
differential operators for all $(E,\nabla)$. By analogy with classical results of
Riemannian geometry [Ep] [St], we identify, in Theorem 5.3.4, sections of $\bbD_X$ with
such natural operators satisfying an extra assumption of regularity. It is extremely interesting
to study natural systems of pseudo-differential operators because intuitively they are given by
kernels which are measures on $\Pi_X$, see (8.1). 

\vskip .3cm

\noindent {\bf (0.4)} The paper is written on two levels. The first level is that of 
elementary differential geometry of bundles, connections and curvatures on $C^\infty$-manifolds.
Our main constructions make sense and seem interesting already at this level and,
with very few exceptions, we stick to this level until \S 6. 

\vskip .2cm

The other level is that of formal geometry and ind-schemes, necessary to construct $\widehat{\Pi}_X$. For
convenience of the reader, most of the technical issues related to ind-schemes and formal groupoids,
were put into the Appendix. It is here that we address the issue of formal integration of Lie
algebroids, i.e., prove the groupoid analog of the fact that any Lie algebra over a field of characteristic
0 gives rise to a formal group. Because integration to a Lie groupoid, even locally, is not always possible, 
see [Mack], we discuss formal integration systematically. Unlike the Lie algebra case,
working directly with the Campbell-Hausdorff series can be confusing here. Our approach,
although essentially equivalent, is based on dualization of $U(\cG)$, the enveloping
algebra of a Lie algebroid $\cG$. The corresponding fundamental structure on $U(\cG)$
is that of a {\em bialgebroid} as defined by J.-H. Lu [Lu] and P. Xu [Xu]. In our situation we need
a  particular case when the cosource and cotarget maps coincide and the algebra of ``objects''
is commutative but not central. This case (but including the antipode which we do not require)
was considered in the recent paper of J. Mr\v cun [Mr].
We call such special bialgebroids {\em left bialgebras} and
give a self-contained exposition for reader's convenience. 

\vskip .3cm

\noindent {\bf (0.5)} The paper is organized as follows. In \S 1 we recall main facts about 
Lie algebroids and their algebraic counterparts, Lie-Rinehart algebras. In \S 2 we present
our main algebraic construction: that of a free Lie-Rinehart algebra generated by an
anchored module. In \S 3 we concentrate on $\cP_X$, the free Lie algebroid generated by
the tangent bundle of a $C^\infty$-manifold $X$. The stabilizer Lie algebra of $\cP_X$ at 
$x\in X$, is the subject of \S 4. In \S 5 we study noncommutative differential operators
and more general natural systems of differential operators in bundles with connections.
In \S 6 we extend the results of \S 5 to the algebro-geometric situation. In \S 7 we introduce
the formal neighborhood $\widehat{\Pi}_X$ and relate it to Kontsevich's moduli spaces.
Finally, \S 8 is devoted to informal discussion of possible further directions
motivated by the present work.

\vskip .3cm

\noindent {\bf (0.6)} I would like to thank J.-H. Lu, J. Mr\v cun,  B. Shoikhet and P. Xu  for pointing out
references to some  previous work relevant to this paper. A large part of this paper was written
during my visit to Max-Planck-Institut f\"ur Mathematik in Bonn in Fall 2006. I would
like to thank the Institute for support and providing excellent working conditions.
In addition, this research was partly supported by the NSF.

\vskip 2cm

\centerline {\bf 1. Reminder on  Lie algebroids.}

\vskip 1cm

\noindent {\bf (1.1) Lie-Rinehart algebras.}
Let $k$ be a field of characteristic 0 and $A$ be a commutative $k$-algebra
with unit. We denote by $\Der(A)$ the $A$-module of derivations of $A$
vanishing on $k$. Recall  that $\Der(A)$ is naturally a Lie $k$-algebra
with respect to the usual commutator. 

 A Lie-Rinehart $A$-algebra is an $A$-module $L$ equipped with
a structure of a $k$-Lie algebra and with a morphism $a: L\to\Der(A)$
of $A$-modules called the {\em anchor map}. 
These data are required to satisfy the following
properties: 

\vskip .2cm

\noindent (1.1.1) $a$ is a morphism of Lie $k$-algebras. 

\vskip .1cm

\noindent (1.1.2) For any  $f\in A$  and sections $x,y\in L$ 
we have
$$[x, fy] - f\cdot [x,y] = a(x)(f)\cdot y.$$

\vskip .2cm

This concept goes back to [Ri]. 
We denote by $\Rin_A$ the category of Lie-Rinehart $A$-algebras.

\vskip .2cm

\noindent {\bf (1.1.3) Examples.} (a) When $A=k$, a Lie-Rinehart $k$-algebra
 is the same
as a Lie $k$-algebra. More generally, a Lie-Rinehart $A$-algebra with the anchor
map being zero is the same as an Lie $A$-algebra in the usual sense. 

\vskip .1cm

(b) $\Der(A)$ with  $a=\id$  is a Lie-Rinehart
$A$-algebra.

\vskip .1cm

(c) Let $M$ be an $A$-module. The {\em Atiyah algebra} of $M$ is the
Lie-Rinehart $A$-algebra $\cA_M$ whose elements are pairs $(\phi, D)$ with
$\phi\in \End_k(M)$ and $D\in\Der(M)$ satisfying the following property:
$$\phi(fm) - f\cdot \phi(m) = D(f)\cdot m, \quad f\in A, m\in M.$$
One sees that commuting the $\phi$'s and the $m$'s makes $\cA_M$
into a $k$-Lie algebra. Further, the anchor map $(\phi, D)\mapsto D$
makes $\cA_M$ into a Lie-Rinehart algebra. Cf. [Kal]. 

\vskip .2cm

For a Lie-Rinehart algebra $L$ we denote by $L^\circ$ the kernel of the
anchor map. This is an Lie $A$-algebra in virtue of (1.1.2).
 A Lie-Rinehart algebra is called transitive, if $a$ is surjective. 

\vskip .1cm

\proclaim (1.1.4) Definition. Let $L$ be a Lie-Rinehart $A$-algebra. An $L$-module
is an $A$-module $M$ together with a morphism of Lie-Rinehart $A$-algebras
$L\to\cA_M$. 

Thus, in particular, $M$ is a module over $L$ as a $k$-Lie algebra.

\vskip .3cm

\noindent {\bf (1.2) Enveloping algebras.} We follow [HS], \S 3.4,
see also [Ri], \S 2.  
Let $L$ be a Lie-Rinehart $A$-algebra. Its (twisted)
enveloping algebra $U(L) = U_A(L)$ is the associative
algebra satisfying the following properties:

\vskip .2cm

\noindent (1.2.1) $L$-modules are the same as left modules over $U_A(L)$ as an associative algebra.

\vskip .1cm

\noindent (1.2.2) $A$ is a subalgebra of $U_A(L)$; there is a natural algebra filtration
$F_\bullet$ of $U_A(L)$ with $F_0 U_A(L)=A$ and
$$\gr^F_\bullet U_A(L) = S^\bullet_A(L),$$
the symmetric algebra of $L$ considered as an $A$-module.

\vskip .2cm

The explicit construction of $U_A(L)$ given in [HS], \S 3.5,  is as follows. One starts from $U_k(L)$,
the usual enveloping algebra of $L$ considered as a Lie $k$-algebra, see [Di] [Re].
Let $U_k^+(L)\subset
U_k(L)$ be the augmentation ideal. We consider it as an associative $k$-algebra without unit. 
 Further, let $U_A(L)^+$ be the quotient algebra  of $U_k(L)^+$ by the relations
$$ x\cdot fy - fx\cdot y = a(x)(f)\cdot y, \quad x,y\in L, f\in A.\leqno (1.2.3)$$
Then we define $U_A(L) = U_A(L)^+ \oplus A$ with the algebra structure given by
$$f\cdot x = fx; 
\quad x\cdot f = fx + a(x)(f), \quad f\in A, x\in L.\leqno (1.2.4)$$
The above properties follow right away from the definition. For example,
$F_i U_A(L) = F_iU_k(L)^+ \oplus A$, where $F_i U_k(L)$ is the canonical filtration
of the usual enveloping algebra. 

\vskip .2cm

\noindent {\bf (1.2.5) Example.} Let $k=\bbR$ and $X$ be a $C^\infty$-manifold. Let $A=C^\infty(X)$
be the algebra of smooth functions on $X$. Then $\Der(A)$ is the Lie algebra of vector fields on $X$,
and $U_A(\Der (A))$ is the algebra of differential operators on $X$.

\vskip .2cm

Note that $A$ is an $L$-module (via the anchor map), and thus a left $U_A(L)$-module. 
We have the augmentation map
$$\epsilon: U_A(L)\to A, \quad P\mapsto P\cdot 1,\leqno (1.2.6)$$
which is a morphism of left $U_A(L)$-modules. Alternatively, $\epsilon$ can be described
as the projection to $A$ along $U_A(L)^+$. 

\vskip .3cm

\noindent {\bf (1.3) Lie algebroids: $C^\infty$-setting.} Let $X$ be a $C^\infty$-manifold.
By a vector bundle on $X$ we mean a locally trivial $\bbR$-vector bundle, possibly of infinite rank.
In other words, a vector bundle $E$ is the same as a locally free sheaf of modules over
$C^\infty_X$, the sheaf of $C^\infty$-functions.

Let $T_X$ be the tangent bundle of $X$. As before, sections of $T_X$ (i.e., vector fields)
form a Lie-Rinehart $C^\infty(X)$-algebra (with $k=\bbR$). So we will call a Lie algebroid on $X$
a vector bundle $\cG$ on $X$ with a morphism of vector bundles $\alpha: \cG\to T_X$
and a structure of a Lie algebra in sections of $\cG$ so that $\alpha$ preserves the
bracket and  the analog of (1.1.2) is satisfied.

 In particular, for a vector bundle $E$ on $X$ we have the Atiyah algebroid $\cA_E$. 
It can be defined more classically as a sheaf of differential operators $P: E\to E$
of order $\leq 1$ whose first order symbol lies in $T_X\otimes 1\subset T_X\otimes\End(E)$.

A module over a Lie algebroid $\cG$ is a vector bundle  $E$, assumed here
 to be  {\em of finite rank}, equipped 
with a morphism of Lie algebroids
$\cG\to \cA_E$. 

\vskip .2cm

\noindent {\bf (1.3.1) Example.} A module over a Lie algebroid $T_X$ is just a vector
bundle on $X$ equipped with a flat connection. 

\vskip .2cm

If $\cG$ is transitive, i.e., if $\alpha$ is surjective, then $\cG^\circ = \Ker(\alpha)$
is a bundle of Lie algebras, so every fiber is of it is an $\bbR$-Lie algebra. 

We also have the concept of the enveloping algebra $U(\cG)$ of a Lie algebroid. Thus
$U(T_X) = \cD_X$ is the sheaf of $C^\infty$-differential operators on $X$.

\vskip 2cm

\centerline {\bf 2. Free Lie-Rinehart algebras and Lie algebroids}

\vskip 1cm

\noindent {\bf (2.1) Free Lie-Rinehart algebras. } 
 One version of the concept of free Lie-Rinehart algebras was introduced in
[CLP]. Here we define a different concept. The relation of our construction
with that of [CLP] will be explained in (2.2.7). 

 Let $k$ be a field of characteristic
zero and $A$ be a commutative $k$-algebra. By an anchored $A$-module we mean
an $A$-module $M$ together with a homomorphism of $A$-modules $b: M\to
 \Der(A)$. Anchored $A$-modules form a category in an obvious way. We denote
this category by $\Anch_A$. Thus we have the forgetful functor
$$\phi: \Rin_A \to \Anch_A,\leqno (2.1.1)$$
which takes a Lie-Rinehart $A$-algebra $L$ into $L$ considered as just
an anchored $A$-module. 

\proclaim (2.1.2) Theorem. The functor $\phi$ has a left adjoint
functor $\FR: \Anch_A\to\Rin_A$ whose value on an anchored $A$-module $M$
is called the free Lie-Rinehart algebra generated by $M$. Thus we have
natural isomorphisms
$$\Hom_{\Rin_A}(\FR(M), L) = \Hom_{\Anch_A}(M, L), \quad M\in\Anch_A, 
L\in\Rin_A.$$

\vskip .2cm

\noindent {\bf (2.1.3) Example.} Suppose that the map $b$ in an anchored
$A$-module $M$ is zero. Then $\FR(M) = \FL(M/A)$ is the free Lie 
$A$-algebra generated by the $A$-module $M$. 

\vskip .2cm

 Let us recall, for the purpose of the future 
generalization,
some details of the construction of $\FL(M/A)$. To construct it,
we start with $\FL(M/k)$, the free Lie $k$-algebra generated by $M$ as
a $k$-vector space. See [Re] for  background on free Lie algebras
over a field. 
 The Lie $A$-algebra $\FL(M/A)$ is obviously
the  quotient of $\FL(M/k)$ by the ``$A$-linearity relations''. 
In order to impose these relations in a systematic way, we consider
the natural grading
$$\FL(M/k) = \bigoplus_{d=1}^\infty \FL_d(M/k), \leqno (2.1.4)$$
where $\FL_d(M/k)$ is the space spanned by brackets involving exactly
$d$ elements of $M$. For example,
$$FL_1(M/k) = M, \quad FL_2(M/k) = \Lambda^2(M).\leqno (2.1.5)$$
Then we construct $FL(M/A)$  as a graded $A$-Lie algebra:
$$\FL(M/A) = \bigoplus_{d=1}^\infty \FL_d(M/A)\leqno (2.1.6)$$
and define the graded components (which are $A$-modules) inductively,
 starting from $\FL_1(M/A)=M$. 

Suppose that for $i=1, ..., d$ we have defined the $A$-modules
$FL_i(M/A)$ and surjections of $k$-vector spaces
$p_i: FL_i(M/k)\to FL_i(M/A)$. 
We define $FL_{d+1}(M/A)$ as the quotient of $FL_{d+1}(M/k)$ be the 
following relations:
$$[x,r], \quad x\in M, r\in\Ker(p_d); \leqno (2.1.7)(a)$$
$$[fx, y] = [x, fy], \quad x\in M, y\in \FL_d(M/A). \leqno (2.1.7)(b)$$
Here the brackets on the left in (2.1.7)(b) are well defined modulo
the relations from (2.1.7)(a). Then we make $\FL_{d+1}(M/A)$ into an
$A$-module by defining $f[x,y]$ to be the common value of the two
brackets in (2.1.7)(b). 

 \vskip .3cm

\noindent {\bf (2.2) Proof of Theorem 2.1.2.}

We generalize the approach to the construction of $\FL(M/A)$ outlined
in Example 2.1.3.  The only difference is that we work with filtered and not graded
objects.

 Thus, we construct $\FR(M)$ as the union of an increasing
sequence of anchored $A$-modules
$$\FR_{\leq 1}(M)\subset \FR_{\leq 2}(M)\subset ... \leqno (2.2.1)$$
which are defined inductively strarting from $\FR_{\leq 1}(M) = M$. 
Suppose we have defined an anchored $A$-module
$\FR_{\leq d}(M)\buildrel a_d\over\longrightarrow \Der (A)$ and
a surjective homomoprhism of $k$-vector spaces 
$$\FL_{\leq d}(M/k) = \bigoplus_{i=1}^d \FL_i(M/k)\buildrel q_d
\over\longrightarrow \FR_{\leq d}(M).$$
We define $\FR_{\leq d+1}(M)$ as the quotient of 
$\FL_{\leq d+1}(M/k)$ by the following relations:
$$[x,r], \quad x\in M, r\in\Ker(q_d); \leqno (2.2.2)(a)$$
$$[fx,y] - [x, fy] = a_d(y)(f)\cdot x - b(x)(f)\cdot y, \quad
x\in M, y\in\FR_{|leq d}(M). \leqno (2.2.2)(b)$$
Here, as before, the brackets on the left in (2.2.2)(b) are well
defined modulo the relations from (2.2.2)(a). 

The $A$-module structure in $\FR_{\leq d+1}(M)$ is defined by
$$f[x,y] = -a_d(y)(f)\cdot x + [fx, y] = 
[x, fy] + b(x)(f)\cdot y, \quad x\in M, y\in\FR_{\leq d}(M).
\leqno (2.2.3)$$
The anchor map is defined by 
$$a_{d+1}([x,y]) = [b(x), a_d(y)], \quad x\in M, y\in\FR_{\leq d}(M).
\leqno (2.2.4)$$
The Lie algebra structure on $\FR(M) = \bigcup_d \FR_{\leq d}(M)$
is induced by that on $\FL(M/k)$. We leave to the reader the verification
that $\FR(M)$ is a Lie-Rinehart algebra and satisfies the
adjunction as in Theorem 2.1.2. \qed

\vskip .2cm

The construction of $\FR(M)$ inplies the following:

\proclaim (2.2.5) Proposition. (a) The filtration $\{\FR_{\leq d}(M)\}$
makes $\FR(M)$ into a filtered Lie-Rinehart algebra, and the associated graded
Lie-Rinehart algebra is isomorphic to $\FL(M/A)$ with trivial anchor map . \hfill\break
(b) Let $\FR^\circ (M)\subset\FR(M)$ be the kernel of the anchor map
(so it is an $A$-Lie algebra). Then the induced filtration on $\FR^\circ(M)$
is compatible with the $A$-Lie algebra structure, and the associated graded
$A$-Lie algebra is isomorphic to $\FL_{\geq 2}(M/A)$:
$$\bigoplus_{d=2}^\infty \FR^\circ_{\leq d}(M)/\FR^\circ_{\leq d-1}(M) = 
\FL_{\geq 2} (M/A), \quad \FR^\circ_{\leq 1}(M)=0.$$

\vskip .2cm

Assume that $A$ is finitely generated over $k$ and $A\to A'$ be an \'etale
extension of $k$-algebras. Then, as well known, $A'$ is finitely generated
too, and
$$\Der(A') = \Der(A)\otimes_A A'.$$            
Because the construction of $\FR(M)$ involves correction terms given by
the action of derivations of $A$, we get the following conclusion.

\proclaim (2.2.6) Proposition. The construction of the free Lie-Rinehart
algebra is compatible with \'etale base change. In other words, if $(M,b)$
is an anchored  $A$-module, and $(M', b')$ its extension to $A'$, with
$M'=M\otimes_A A'$, then
$$\FR(M') = \FR(M)\otimes_A A'.$$
In particular, the construction of $\FR(M)$ is compatible with 
localization.

\vskip .2cm

\noindent {{\bf (2.2.7) Relation to [CLP].} Casas, Ladra and Pirashvili define, in
 {\em loc. cit.} a different concept of free Lie-Rinehart algebras. They start
not with anchored $A$-modules but with {\em $A$-anchored $k$-vector spaces}. By definition,
an $A$-anchored $k$-vector space is a vector space $V$ together with a $k$-linear map
$c: V\to \Der(A)$. Such objects for a category $\Anch_{k/A}$ and we obviously have
the forgetful functor
$$\psi: \Rin_A\to \Anch_{k/A}.$$
Their free Lie-Riehart algebra functor is the left adjoint functor to $\psi$. Let us denote it
by $\FR_{CLP}$. It is clear that the functor $\psi$ factor through the forgetful functor
$\phi$ to the category of anchored $A$-modules. Further, if $(V,c)$ is an $A$-anchored
$k$-vector space, then we have an anchored $A$-module $(V\otimes_k A, c\otimes 1)$
freely generated by $V$, and we have:

\proclaim (2.2.8) Proposition. If $(V,c)$ is an $A$-anchored $k$-vector space, then
$$\FR_{CLP}(V) = \FR(V\otimes_k A).$$

Our construction has an additional flexibility in that it is defined for not necessarily free
$A$-modules and localizes on $\Spec(A)$ if $A$ is finitely generated.

\vskip .3cm

\noindent {\bf (2.3) Free Lie algebroids: $C^\infty$-setting.}
Let $k=\bbR$. Let $X$ be a $C^\infty$-manifold and $(V, \beta)$ be an anchored
vector bundle on $X$, i.e., a vector bundle (possibly infinite-dimensional)
 with a morphism of vector bundles $\beta: E\to T_X$. Then we get a Lie
algebroid $\cF(V)$, called the free Lie algebroid generated by $E$ with
properties similar to the ones listed above. 

\proclaim (2.3.1) Definition.  Let $X$ be a $C^\infty$-manifold, and
$\beta: V\to T_X$ be an anchored vector bundle. A pre-module over $V$
is a vector bundle $E$ together with a morphism of anchored bundles
$V\to\cA_E$ from $V$ to the Atiyah algebroid of $E$.

\vskip .2cm

\noindent {\bf (2.3.2) Example.} A pre-module over $T_X$ is the same as a 
vector bundle with connection (not necessarily flat).

\vskip .2cm

The universal property of free Lie algebroids implies the following:

\proclaim (2.3.3) Proposition. Modules over $\cF(V)$ are the same as
pre-modules over $V$. 

\vskip 2cm

\centerline {\bf 3. The path Lie algebroid.}

\vskip 1cm

\noindent {\bf (3.1) Main properties.} Let $X$ be a $C^\infty$-manifold.
The free Lie algebroid $\cF(T_X)$ generated by $T_X$ will be called the
{\em path algebroid} of $X$ and denoted by $\cP_X$. The kernel of the
anchor map of $\cP_X$ will be denoted by $\cP^\circ_X$. This is a bundle
of Lie algebras on $X$. Let us summarize the properties of
$\cP_X$ that follow from the general results of \S 2. 

\proclaim (3.1.1) Proposition. (a) Modules over $\cP_X$ are the same
as finite vector bundles on $X$ with connection. \hfill\break
(b) $\cP_X$ is the union of subbundles $\cP_{X, \leq d}$, $d\geq 1$
of finite rank. The filtration $\{\cP_{X, \leq d}\}$ is admissible (A.5.7), and
$$\gr (\cP_X) = \FL(T_X)$$
is the bundle of fiberwise free Lie algebras of $T_X$. 
\hfill\break
(c) The induced filtration $\cP^\circ_{X, \leq d} = \cP_{X, \leq d}\cap
\cP^\circ_X$ has $\cP^\circ_{X, \leq 1}=0$. It is compatible with the
Lie algebra structure in the fibers, and
$$\gr(\cP^\circ_X) = \FL_{\geq 2}(T_X).$$
(d) We have a canonical embedding
$$c: \Lambda^2 (T_X)\isomoto \cP^\circ_{X, \leq 2}\subset \cP_{X, \leq 2},$$
which takes the wedge product $v\wedge w$ of two vector fields into
$$c(v\wedge w) = [v,w]_{\cP} - [v,w]_{\Lie}.$$
Here $[v,w]_{\cP}$ is the bracket in $\cP_X$, while
$[v,w]_{\Lie}$ is the standard Lie bracket of $v,w$ as vector fields.
\qed

\vskip .3cm

\noindent {\bf (3.2) Noncommutative vector fields.} 
Let $\Bun_\nabla(X)$ be the category of vector bundles of finite rank on $X$ equipped
with connections. It is a tensor category with respect to the usual tensor
product $\otimes$. 

\vskip .2cm

\proclaim (3.2.1) Definition. A noncommutative vector field on $X$
is a rule $D$ which to each $(E, \nabla)\in \Bun_\nabla (X)$
associates a differential operator $D_{E, \nabla}: E\to E$ satisfying the
following properties:\hfill\break
(a) Naturality: for any morphism $\phi: (E, \nabla)\to (E', \nabla')$
of bundles with connections we have $D_{E', \nabla'}\circ\phi = \phi\circ D_{E,\nabla}$.
\hfill\break 
(b) Leibniz rule with respect to the tensor product: 
for any two bundles with connections $(E,\nabla)$ and $(E', \nabla')$
we have
$$D_{E\otimes E', \nabla\otimes \nabla'} = D_{E,\nabla}\otimes 1 +
 1\otimes D_{E', \nabla'}.$$

Clearly, noncommutative vector fields form a sheaf on $X$ which
we denote by $\VVect_X$. Further, $\VVect_X$ is a sheaf of Lie algebras
with respect to the usual commutator of differential operators in
any $(E,\nabla)$. 

If $D$ is an noncommutative vector field, then taking $(E,\nabla) = (E', \nabla')=
 C^\infty_X$
to be the trivial rank 1 bundle with standard connection, we find
from the $\otimes$-Leibniz rule that
$$D_{C^\infty_X}: C^\infty_X\to C^\infty_X$$
is a ring derivation, i.e., it is a vector
field in the usual sense. This defines a morphism of sheaves of Lie
algebras $\alpha: \VVect_X\to T_X$. 

\proclaim (3.2.2) Proposition. The map $\alpha$ makes $\VVect_X$ into
a Lie algebroid (sheaf of Lie-Rinehart algebras). \qed

Further, applying the $\otimes$-Leibniz rule to $C^\infty_X\otimes E = E$, we find
that each $D_{E, \nabla}: E\to E$ is a first order differential operator and its first
order symbol is $\alpha(D)\otimes 1_E$. In other words:

\proclaim (3.2.3) Proposition. For each $(E, \nabla)\in\Bun_\nabla (X)$ the correspondence
$D\mapsto D_{E, \nabla}$ defines a morphism of Lie algebroids
$$\VVect_X\to \cA_E.$$

\vskip .2cm

\noindent {\bf (3.2.4) Examples.} (a) Let $v$ be a usual vector field
on $X$. Then the rule $D_v$ which to each $(E,\nabla)$ associates
$\nabla_v: E\to E$ (covariant derivative along $v$), is a noncommutative
vector field. 

\vskip .1cm

(b) Let $\xi$ be a bivector field on $X$, i.e., a section of
 $\Lambda^2(T_X)$. Then we have the noncommutative vector field $\Phi_\xi$
which to each $(E,\nabla)$ associates the endomorphism
$(F_\nabla, \xi): E\to E$. Here $F_\nabla\in\Omega^2_X\otimes\End(E)$
is the curvature of $\nabla$. 

\vskip .1cm

(c) Note that the correspondence $v\mapsto D_v$ from (a) is not a Lie algebra
homomorphism. Indeed, we have
$$[D_v, D_w] = D_{[v,w]} + \Phi_{v\wedge w}.$$
This is just the definition of the curvature. 

\proclaim (3.2.5) Theorem. The correspondence $v\mapsto D_v$ extends
to a monomorphism of Lie algebroids
$$h: \cP_X\to \VVect_X.$$

The existence of $h$ follows from the universal property of $\cP_X$, as
we have  a morphism of anchored vector bundles $T_X\to \VVect_X$. 
The proof that $h$ is an monomorphism will be given in \S 5 below where
we will also characterize the image. .

\vskip .3cm

\noindent {\bf (3.3) Example: the case of $\bbR^n$.}
  Let $X=\bbR^n$ with coordinates $x_1, ...,
x_n$. Let us describe $\cP(\bbR^n)$, the space of global sections of $\cP_{\bbR^n}$
explicitly. 
Let $\partial_i=\partial/\partial x_i$ and $D_i$ be the section
of $\cP_{\bbR^n} = \cF(T_{\bbR^n})$ corresponding to the section
$\partial_i$ of $T_{\bbR^n}$. 

Consider the associative algebra $\bbD(\bbR^n)$ generated by 
 $C^\infty$-functions $f(x_1, ..., x_n)$
and the symbols $D_i$ subject only to the relations
$$[D_i, f] = {\partial f\over\partial x_i}. \leqno (3.3.1)$$
Note that we do not assume the $D_i$ to commute with each other; in
fact, they generate the free associative algebra $\bbR\langle D_1, ..., D_n\rangle$. 

Let $\bbD(\bbR^n)_{\Lie}$ be $\bbD(\bbR^n)$ considered as a Lie algebra with
respect to the usual commutator $[x,y] = xy-yx$. Then the Lie subalgebra in
$\bbD(\bbR^n)_{\Lie}$ generated by $D_1, ..., D_n$ is the free Lie algebra
$\FL(\bbR^n)$. 

The space $\cP(\bbR^n)$ is identified with 
the Lie subalgebra in $\bbD(\bbR^n)_{\Lie}$ generated
by elements of the form
$$\sum_{i=1}^n f_i(x_1, ..., x_n) D_i.$$
Indeed, such an element corresponds to $D_v$ from Example 3.2.4(a), where
$v$ is the vector field $\sum f_i \partial_i$. 

We see easily that a general element of $\cP(\bbR^n)$ 
 can be uniquely written in the form (sum with only finitely many
nonzero summands)
$$\eta = \sum_{p=1}^\infty \sum_{i_1, ..., i_p} f_{i_1, ..., i_p}(x) [D_{i_1},
[D_{i_2}, ... [D_{i_{p-1}}, D_{i_p}]...],\leqno (3.3.2) $$
where, for each $p$, the  $(i_1, ..., i_p)$ run over the set of multiindices such that
the brackets in the above sum form a  Hall basis in $\FL_p(\bbR^n)$, 
see [Re]. The structure of a vector bundle, i.e., of a module over
the ring $C^\infty(\bbR^n)$, is given by multiplying the
$f_{i_1, ..., i_p}$ simultaneously with a given function $f$. The
anchor map takes a sum above into $\sum f_i(x) \partial_i$, with
$f_i$ coming from the sub-sum with $p=1$.

Note that the higher commutators of the $D_i$ commute with functions in virtue of (4.3.1).
Let us denote by 
$$\cP_p(\bbR^n)  = C^\infty(\bbR^n)\otimes_{\bbR}
FL_p(\bbR^n)\leqno (3.3.3)$$ 
the part of $\cP(\bbR^n)$ corresponding to the
summand in (3.3.2) with given $p$. Then we have that the degree $\geq 2$ part forms a graded
Lie algebra:
$$[\cP_p(\bbR^n), \cP_q(\bbR^n)] \subset \cP_{p+q}(\bbR^n), \quad p,q \geq 2,\leqno (3.3.4)$$
while
$$[\cP_1(\bbR^n), \cP_q(\bbR^n)] \subset \cP_q(\bbR^n) \oplus \cP_{1+q}(\bbR^n). \leqno (3.3.5)$$

\vskip .1cm

If $E$ is  a trivial rank $N$ bundle on $\bbR^n$ with connection $\nabla$, we have the
connection operators
$$\nabla_i = \partial_i + A_i(x), \quad A_i(x)\in C^\infty(\bbR^n)\otimes \Mat_N(\bbR).
\leqno (3.3.6)$$
The element $D_i\in \cP(\bbR^n)$ is sent by the homomorphism $h$ from the formulation of Theorem
3.2.5, to the noncommutative vector field $D_{\partial_i}$ which acts on $(E,\nabla)$ by
$$(D_{\partial_i})_{E,\nabla} = \nabla_i.\leqno (3.3.7)$$

\vskip .2cm

\noindent {\bf (3.3.8) Remarks.} (a)  The above description of $\cP(\bbR^n)$ is  a particular case of the
construction in [CLP], Example 2.2 called the transformational  Lie-Rinehart algebra. 
This construction was used in {\em loc. cit.} to define the  free Lie-Rinehart algebra
generated by an anchored vector space (the functor $\FR_{CLP}$, see Remark 2.2.7). In fact,
in our case $T_{\bbR^n}$ is trivialized, so $\cP(\bbR^n) = \FR_{CLP}(\bbR^n)$
where $\bbR^n$ on the right is considered as a $C^\infty(\bbR^n)$-anchored vector space by
 sending the
$i$th basis vector of $\bbR^n$ to $\partial_i$. 

\vskip .1cm

(b) The algebra $\bbD(\bbR^n)$ is  identified with 
the enveloping algebra of the
Lie-Rinehart algebra $\cP(\bbR^n)$.  It is
a particular case of algebras of noncommutative differential operators
to be studied later.

\vskip 2cm

\centerline {\bf 4. The fundamental Lie algebra.}

\vskip 1cm

\noindent {\bf (4.1) Interpretation via curvatures.}
Let $X$ be a $C^\infty$-manifold and $x\in X$. We call the {\em fundamental Lie algebra}
of $X$ at $x$ and denote $\cP(X,x)$ 
 the fiber of $\cP^\circ_X$
at $x$. It can be seen as the Lie algebra of  the group
of $x$-based loops in $X$ modulo reparametrizations and cancellations. Proposition
3.1.1 implies the following.

\proclaim (4.1.1) Proposition. The Lie algebra $\cP(X,x)$ comes with a natural filtration
$\{ \cP_{\leq d}(X,x)\}$, $d\geq 2$,  compatible with the Lie algebra structure. The lowest
term of this filtration, i.e., $\cP_{\leq 2}(X,x)$,  is canonically identified with $\Lambda^2(T_xX)$.
Further, the associated graded Lie algebra is canonically identified with the degree $\geq 2$
part of the free Lie algebra generated by $T_xX$:
$$\bigoplus_{d\geq 2} \cP_{\leq d}(X,x) /\cP_{\leq d-1}(X,x)\quad  \simeq \quad \FL_{\geq 2}(T_xX).$$

\vskip .1cm

Another interpretation of $\cP(X,x)$ is that it is the universal receptacle of all the
covariant derivatives of the curvatures of all connections in vector bundles on $X$. 
Indeed, if $(E,\nabla)\in\Bun_\nabla (X)$, then we have the ``holonomy morphism''
of Lie algebras
$$H_{E,\nabla}: \cP(X,x) \to \End(E_x). \leqno (4.1.2)$$
Proposition 4.1.1 can be interpreted as follows. The curvature $F=F_\nabla$ of a connection
$\nabla$ is a tensor $F\in\Omega^2_X\otimes \End(E)$, and the lowest term of the
filtration on $\cP(X,x)$ is mapped via $F$:
$$H_{E,\nabla}|_{\Lambda^2 T_xX} = F|_x: \Lambda^2 T_xX\to \End(E).\leqno (4.1.3)$$
Higher covariant derivatives of $F$ are not invariantly defined, as the bundle $\Omega^2_X\otimes\End(E)$
carries no natural connection induced by $\nabla$. What is naturally defined, is the ``de Rham differential''
$$\nabla^{(2)}: \Omega^2_X\otimes \End(E)\to\Omega^3_X\otimes\End (E),\leqno (4.1.4)$$
which satisfied $\nabla^{(2)}(F_\nabla)=0$, the Bianchi identity. If we choose local coordinates
$x_1, ..., x_n$, we can trivialize both $\Omega^1_X$ and $\Omega^2_X$ and get a map
$$\widetilde{\nabla}^{(2)}:  \Omega^2_X\otimes\End(E)\to \Omega^1_X\otimes\Omega^2_X\otimes\End(E),\leqno (4.1.5)$$
of which $\nabla^{(2)}$ is obtained by antisymmetrization. The image $\widetilde{\nabla}^{(2)}(F_\nabla)$ lies then in
$$\Ker\biggl\{ \Omega^1_X\otimes\Omega^2_X\otimes \End(E) \to \Omega^3_X\otimes\End(E)\biggr\} = 
\Hom\bigl( \FL_3(TX), \End(E)).\leqno (4.1.6)$$
The last identification corresponds to the well known intepretation of the Bianchi identity as the
Jacobi identity for the connection operators $\nabla_i$. However, the operator
$\widetilde{\nabla}^{(2)}$ is not invariantly defined,
and what we have instead is Proposition 4.1.1. 

\vskip .2cm
Similarly to Theorem 3.2.5, elements of $\cP(X,x)$ can be viewed as derivations (in the sense of [KM], 
\S 1.7) of the fiber functor $(E,\nabla) \mapsto E_x$ from $\Bun_\nabla(X)$ to vector spaces.

\vskip .3cm

\noindent {\bf (4.2) Homology of $\cP(X,x)$: preliminaries.} As $\cP(X,x)$ is isomorphic 
(non-canonically)
to a subalgebra of a free Lie algebra, namely, to  $\FL_{\geq 2}(T_xX)\subset \FL(T_xX)$,
it is itself free by the Shirshov-Witt theorem, see [Re], Ch.2. The number of free generators
of $\cP(X,x)$ is, however, infinite. 

Let $k$ be a field of characteristic 0. 
For an Lie $k$-algebra $\Gg$ we denote by $H_\bullet^{\Lie}(\Gg)$ the Lie algebra homology of $\Gg$
with trivial coefficients $k$. Because $\cP(X,x)$ is free (with $k=\bbR$),
we have
$$H_p^{\Lie} \cP(X,x)=0, \quad p\geq 2.\leqno (4.2.1)$$
Let us determine the space
$$H_1^{\Lie} \cP(X,x) = \cP(X,x)/[\cP(X,x), \cP(X,x)],$$
i.e., the maximal abelian quotient. Its importance stems from the following
interpretation of $H_1^{\Lie}$ as ``the space of generators''.

\proclaim (4.2.2) Proposition. Let $\Gg$ be an $\bbZ_+$-graded Lie $k$-algebra which is free
(i.e., has a system of homogeneous free generators). Suppose
$W\subset \Gg$ is a graded  $k$-vector subspace such that the composite map
$$W\to\Gg \to H_1^{\Lie}(\Gg)$$
is an isomorphism. Then $W$ is a space of free generators for $\Gg$, i.e., the natural morphism of
Lie algebras $\FL(W)\to\Gg$ is an isomoprhism. \qed.

\vskip .2cm

We start with recalling the description
of  $H_1^{\Lie}(\FL_{\geq 2}(V))$ where $V$ is a finite-dimensional
$k$-vector space. This homology space is graded:
$$H_1^{\Lie}(\FL_{\geq 2}(V)) = \bigoplus_{d=2}^\infty H_{1,d}^{\Lie}(FL_{\geq 2}(V)),\leqno (4.2.3)$$
where $H_{1,d}^{\Lie}$ is the image of $\FL_d(V)\subset \FL_{\geq 2}(V)$
in the maximal abelian quotient of $FL_{\geq 2}(V)$.

\vskip .2cm

Let  $\Vect_k$ be the category of finite-dimensional
$k$-vector spaces. For any sequence of integers $\alpha = (\alpha_1 \geq ...\geq \alpha_n\geq 0)$
(with arbitrary $n$) we have the polynomial functor $\Sigma^\alpha: \Vect_k\to\Vect_k$, see
[Macd]. For $V=k^n$ the space $\Sigma^\alpha(V)$ is the space of irreducible representation of 
the algebraic group $GL_n/k$ with highest weight $\alpha$.

\proclaim (4.2.4) Theorem. (a)  We have an identification of functors of $V$:
$$H_{1,d}^{\Lie} (\FL_{\geq 2}(V)) \simeq \Sigma^{d-1,1}(V).$$
(b) Let $z_1, ..., z_n$ be a basis of $V$.  Then the images of the  elements
$$[z_{i_1}, [z_{i_2}, ... [z_{i_{d-1}}, z_{i_d}]...], \quad i_1\geq i_2 \geq ... \geq i_{d-1} < i_d,$$
form a  basis of $H_{1,d}^{\Lie}(\FL_{\geq 2}(V))$. \hfill\break
(c) The elements as above, taken for all $d$, form a system of free generators of
$\FL_{\geq 2}(V)$.

Parts (a) and (b) are due to Reutenauer  [Re], \S\S 8.6.12 and 5.3.
Part (c)  follows from  Proposition 4.2.2. \qed

\vskip .2cm

Let us now give a geometric interpretation of Theorem 4.2.4. 
 Consider $V$ as an affine algebraic variety over $k$
and let $\Omega^p(V)$ be the space of regular $p$-forms on $V$. Thus
$$\Omega^p(V) = \bigoplus_{d=0}^\infty \Lambda^p(V^*)\otimes S^d(V^*).$$
The de Rham differential $d: \Omega^p(V)\to\Omega^{p+1}(V)$ makes $\Omega^\bullet (V)$
into a complex exact everywhere except the $0$th term, and
the space of closed $p$-forms has the following well known decomposition as a $GL(V)$-module:
$$\Omega^{p, cl}(V) = d(\Omega^{p-1}(V)) = \bigoplus_{d=1}^\infty \Sigma^{d, 1, ..., 1}(V^*)
\leqno (4.2.5)$$
(with $(p-1)$ occurrences of 1 in the RHS). See, e.g., [GKZ], Proposition 14.2.2. 
 Thus
$$H_1^{\Lie}(\FL_{\geq 2}(V^*)) = \Omega^{2, cl}(V). \leqno (4.2.6)$$
The first homology of $\FL_{\geq 2}(V)$ is the restricted dual of this space with respect to the
decomposition (4.2.5), i.e., the direct sum of the duals of the summands.
Let us denote by
$$\Gamma_p(V) = \bigoplus_{d=0}^\infty \Lambda^p(V)\otimes S^d(V)\leqno (4.2.7)$$
the restricted dual of $\Omega^p(V)$. The de Rham differential $d$ on $\Omega^\bullet(V)$
induces differentials $\partial: \Gamma_p(V)\to\Gamma_{p-1}(V)$ by dualization. The dual of
$$\Omega^{2, cl}(V) = \Ker\bigl\{ d: \Omega^2(V)\to\Omega^3(V)\bigr\}$$
is then
$$\Coker\bigl\{ \partial: \Gamma_3(V)\to\Gamma_2(V)\bigr\} \quad = \quad 
\Ker \bigl\{ \partial: \Gamma_1(V)\to\Gamma_0(V)\bigr\}, \leqno (4.2.8)$$
which we denote by $\Gamma_1^{cl}(V)$. 

\vskip .3cm

\noindent {\bf (4.3) Homology of $\cP(X,x)$ and currents.} 
Let now $k=\bbR$, let $X$ be a $C^\infty$-manifold and $x \in X$. Denote by $\Gamma_p(X,x)$
the space of $p$-currents on $X$ is the sense of de Rham, which are supported at $x$.
In other words, an element of $\Gamma_p(X,x)$ is a linear functional $\Omega^p(X)\to\bbR$
on the space of all $C^\infty$ forms, which depends only on a finite jet of a form at $x$.
It is clear that for an $\bbR$-vector space $V$ considered as a $C^\infty$-manifold the space 
$\Gamma_p(V, 0)$ is the same as $\Gamma_p(V)$ defined earlier. 

As before, the de Rham differential on forms defines differentials
$\partial: \Gamma_p(X,x)\to \Gamma_{p-1}(X,x)$. 

\vskip .1cm

\proclaim (4.3.1) Theorem. The space $H_1^{\Lie}(\cP(X,x))$ is canonically identified with
$\Gamma^{cl}_1(X, x)$, the space of closed 1-currents supported at $x$. 

By dualizing, we get an essentially equivalent formulation:

\proclaim (4.3.2) Corollary. The space $H^1_{\Lie}(\cP(X,x))$ of first cohomology with coefficients
in $\bbR$ is canonically identified with $\widehat{\Omega}^{2, cl}_{X,x}$, the space of
formal germs of closed 2-forms on $X$ near $x$. 

\vskip .2cm

To prove Theorem 4.3.1, we first construct a natural linear map
$$\tau: \cP(X,x)\to\Gamma_1^{cl}(X,x),$$
vanishing on the commutators. Indeed, an element $\zeta\in\Gamma_1^{cl}(X,x)$ can be seen, as above,
as a functional on $\Omega^{2, cl}(X)$ whose value on a closed form $\omega$ depends only on
a finite jet of $\omega$ near $x$. On the other hand, any $\eta\in\cP(X,x)$
gives rise to the holonomy element
$$H_{E,\nabla}(\eta)\in\End(E_x)$$
for any bundle with connection $(E,\nabla)$ on $X$. Let us restrict to the case when the rank of $E$
is equal to 1. Then the curvature $F_\nabla$ is a scalar closed 2-form: $F_\nabla\in\Omega^{2, cl}(X)$.
At the same time $\End(E_x)=\bbR$. 

So by sending $\eta$ to $H_{E,\nabla}(\eta)$ where $(E,\nabla)$ is a line bundle with
connection such that $F_\nabla=\omega$, we associate to $\eta$ a linear functonal
of $\omega\in\Omega^{2, cl}(X)$ which clearly depends only on a finite jet of $\omega$ at $x$.
This defines the map $\tau$.

Next, to show that $\tau$ is an isomorphism, it is enough to choose a local coordinate
system
on $X$ near $x$, in which case we can assume $X=V=\bbR^n$ and the homomoprhism $\tau$ is identified with
that of Theorem 4.2.4. \qed

\vskip .2cm

\noindent {\bf (4.3.3) Remarks.}(a)  Theorem 4.3.1 has the following group-theoretical analog.
Let $F_n$ be the free group on $n$ generators, and $[F_n, F_n]$ be its commutator. Being
a subgroup of a free group, it is free itself and has only the first group homology space
$H_1^{\operatorname{Gr}}([F_n, F_n], \bbR)$ nontrivial. This space is described as follows. Recall that $BF_n$,
the classifying space of $F_n$, is the bouquet of $n$ circles. Let $\bbR^n_\square$ be the
Euclidean space $\bbR^n$ with CW-decomposition into unit cubes of the standard integer lattice
$\bbZ^n$. Then $B([F_n, F_n]) = {\operatorname{Sk}}_1(\bbR^n_\square)$ is the 1-skeleton of
this CW-complex.  Accordingly,
$$H_1^{\operatorname{Gr}}([F_n, F_n], \bbR) = H_1^{\operatorname{Top}}({\operatorname{Sk}}_1(\bbR^n_\square),
\bbR) = Z_1(\bbR^n_\square, \bbR)$$
is the space of cellular 1-cycles of the CW-complex $\bbR^n_\square$. 

\vskip .1cm

(b) The identification in (4.2.6) and thus in Theorem 4.3.1, seems  related to the results of
Feigin and Shoikhet [FS] on the structure of $[A,A]/[A, [A,A]]$ where $A$ is the
free associative algebra generated by  $V$. 

\vskip .1cm

(c) An identification equivalent to that of
  Theorem 4.2.4 (a) was also found by Polyakov ([Po], \S 3)  who was
studying ``gauge invariant words'', i.e., traces of iterated covariant
derivatives of the curvature of an indeterminate connection on $X$. Such iterated
derivatives are labelled, as discussed in (4.1), by elements of $\cP(X,x)$,
while traces vanish on $[\cP(X,x), \cP(X,x)]$.

\vskip .3cm

\noindent {\bf (4.4) Formal classification of connections.} 
Let $k$ be a field of characteristic $0$, and
$\GD_n = \Spf \, k[[x_1, ..., x_n]]$
be the $n$-dimensional formal disk over $k$, see Example A.2.2(a) . We can then speak about (finite rank)
vector bundles with connection $(E,\nabla)$ on $\GD_n$. Every vector bundle on $\GD_n$ is trivial
but perhaps not canonically.
If $E = \cO^N$
is the standard trivial rank $N$ bundle, a connection in $E$ is given by operators
$$\nabla_i= \partial_i + A_i, \quad A_i \in\Mat_N\, k[[x_1, ..., x_n]],\leqno (4.4.1)$$
and isomorphisms of such connections are given by gauge transformations
$$A_i\mapsto g^{-1} A_i g + g^{-1} \partial_i g, \quad g\in GL_N\, k[[x_1, ..., x_n]].
\leqno (4.4.2)$$
We denote by $\Bun_\nabla (\GD_n)$ the category of finite rank vector bundles with connections on $\GD_n$. 

\vskip .2cm

More generally, let $G$ be an affine algebraic group over $k$ with Lie algebra $\Gg$. 
Then we can speak about principal $G$-bundles with connection on $\GD_n$. For example, gauge
transformations of the standard trivial $G$-bundle form the group $G(k[[x_1, ..., x_n]])$. 

\vskip .1cm

\proclaim (4.4.3) Theorem. The category $\Bun_\nabla(\GD_n)$ is equivalent
(as a tensor category)  to the category of
finite-dimensional representations of the Lie algebra $\FL_{\geq 2}(k^n)$.

This is a corollary of Theorem A.7.3 about transitive Lie algebroids on $\GD_n$. \qed

\vskip .1cm

\proclaim (4.4.4) Corollary. Let $G$ be an algebraic group as above. Then gauge equivalence
classes of connections in principal $G$-bundles on $\GD_n$ are in bijection with
$G$-conjugacy classes of homomorphisms $\FL_{\geq 2}(k^n) \to \Gg$. \qed

\vskip .1cm

Let us consider the case of the trivial $G$-bundle $G\times\GD_n$, so a connection
$\nabla$ is given by $\nabla_i=\partial_i + A_i$  with $A_i\in\Gg [[x_1, ..., x_n]]$. By a
{\em restricted gauge transformation} of $G\times\GD_n$ we mean a transformation
$$g(x_1, ..., x_n)\in G(k[[x_1, ..., x_n]])$$
whose value at $0$ (i.e., the constant term) is equal to 1. 
Let
$$F_{ij} = \partial_i A_j -\partial_j A_i - [A_i, A_j]$$
be the $(i,j)$th component of the curvature of $A$.

\proclaim (4.4.5) Corollary (Taylor theorem for connections). 
Given a connection $\nabla$ as above, the set of elements
$$\nabla_{i_1}\nabla_{i_2} ... \nabla_{i_{d-2}} F_{i_{d-1}, i_d} \in \Gg, \quad i_1\geq i_2\geq ... \geq
i_{d-1} < i_d, \,\,\, d\geq 2,$$
determine $\nabla$ uniquely up to restricted gauge equivalence. 
Conversely, given any system of elements
$$\gamma_{i_1, ... , i_d}\in\Gg, \quad  \quad i_1\geq i_2\geq ... \geq
i_{d-1} < i_d, \,\,\, d\geq 2,$$
there exists a unique (up to restricted gauge equivalence) connection $\nabla$ on $\GD_n$
with 
$$\nabla_{i_1}\nabla_{i_2} ... \nabla_{i_{d-2}} F_{i_{d-1}, i_d} = 
\gamma_{i_1, ... , i_d}.$$ 

This is a consequence of Corollary 4.4.4 as well as Theorem 4.3.4(c)  describing a system
of free generators of $FL_{\geq 2}(k^n)$. \qed

\vskip 2cm

\centerline {\bf 5. Noncommutative differential operators.}

\vskip 1cm

\noindent {\bf (5.1) Basic definitions.} Let $X$ be a $C^\infty$-manifold. We denote by
$\bbD_X$ the sheaf $U(\cP_X)$ of enveloping algebras of the Lie algebroid $\cP_X$. Sections
of $\bbD_X$ will be called {\em noncommutative differential operators}. 

By Theorem A.3.11, $\bbD_X$ is a sheaf of left $C^\infty_X$-bialgebras; in particular, we have a morphism
of algebras
$$\Delta: \bbD_X\to \bbD_X\overline{\otimes}_{C^\infty_X} \bbD_X. \leqno (5.1.1)$$
Further,  each section $P$ of $\bbD_X$ gives a differential operator 
$P_{E,\nabla}: E\to E$ for any $(E,\nabla)\in\Bun_\nabla(X)$. These operators are
natural with respect to morphisms of bundles with connections and behave with respect
to the tensor product as follows:
$$P_{E\otimes E', \nabla\otimes\nabla'} = \sum (P^{(1)}_i)_{E,\nabla} \otimes (P^{(2)}_i)_{E', \nabla'}, \quad 
\Delta(P) = \sum P^{(1)}_i\otimes P^{(2)}_i. \leqno (5.1.2)$$
By a {\em natural differential operator} on $X$ of order $\leq d$ we will mean
a   system $\{P_{E,\nabla}: E\to E\}$
of differential operators of order $\leq d$  defined for all bundles with connections
$(E, \nabla)$ at once and satisfying the naturality property as above.
Compare [St], Def. 3.4 for the case of Riemannian manifolds, not connections in bundles
on a given manifold.
 We denote by $\cN\cD_X^{\leq d}$
the sheaf of natural differential operators on $X$ of order $\leq d$, and
by $\cN\cD_X$ the union of these sheaves for all $d$.

By the above, each section of $\bbD_X$ gives a natural differential operator, so
we have a homomorphism of sheaves of algebras on $X$ extending (3.2.5): 
$$h: \bbD_X\to \cN\cD_X. \leqno (5.1.3)$$

\vskip .2cm

\noindent {\bf (5.1.4) Example.}  Let $X=\bbR^n$. Then the algebra $\bbD(\bbR^n)$ of global
sections of $\bbD_{\bbR^n}$ was already described in (3.3). It is generated by functions
$f(x_1,..., x_n)$ and symbols $D_1, ..., D_n$ subject only to the relations (3.3.1). 
For example, we have the (noncommutative) Laplace operator $\Delta\in\bbD(\bbR^n)$ whose
definition and the action on vector bundles with connections are given by:
$$\Delta = \sum_{i=1}^n 
D_i^2, \quad  \Delta_{E, \nabla} = \sum_{i=1}^n \nabla_i^2: E\to E.$$
Further, we have a similar situation for any Riemannian manifold $X$: the metric
defines an element $\Delta_X\in\bbD(X)$.

\vskip .3cm

\noindent {\bf (5.2) Two filtrations on $\bbD_X$.} As for the enveloping algebra of any Lie
algebroid, $\bbD_X=U(\cP_X)$ has the filtration $F$, see  (1.2.2) by the number of noncommutative
vector fields needed to produce a given section $P$ of $\bbD_X$.  Thus
$$\gr^{F}(\bbD_X) = S^\bullet(\cP_X)\leqno (5.2.1)$$
is a commutative algebra. 

\vskip .1cm

Second, recall that $\cP_X = \bigcup_{d=1}^\infty \cP_{\leq d}$ is itself a filtered Lie algebroid, so we can take the
filtration$\{U_{\leq d}(\cP_X)\}$ on $U(\cP_X)$ induced by this filtration on $\cP_X$, as in (A.5.8). We call this
filtration the {\em filtration by degree} and write $\bbD^{\leq d}_X$ for $U_{\leq d}(\cP_X)$. We also write
 $\deg(P)\leq d$, if $P\in \bbD^{\leq d}_X$. 
We then see easily that
$$\gr^{\deg}(\bbD_X) \quad = \quad T^\bullet(T_X)\quad=\quad \bigoplus_{d=0}^\infty
T_X^{\otimes d}\quad  =\quad  U(\FL(T_X)), \leqno (5.2.2)$$
is the full tensor algebra of $T_X$, with the grading induced by the grading on $\FL(T_X)$
on its enveloping algebra.
 So it is noncommutative.     

\vskip .1cm

For a noncommutative differential operator $P$ on $X$ of degree $\leq d$ we denote by
$\Smbl_d(P) \in T_X^{\otimes d}$ the highest symbol of $P$, i.e., the image of $P$ in
$\gr_d^{\deg}(\bbD_X)$. Thus $\Smbl_d(P)$ can be seen as a noncommutative
homogeneous polynomial on $T^*X$, the cotangent bundle of $X$.

\proclaim (5.2.3) Proposition. Let $(E,\nabla)\in \Bun_\nabla(X)$, and $P$ be a section of
$\bbD_X$ with $\deg(P)\leq d$. Then the linear differential operator 
$P_{E,\nabla}:E\to E$ has order $\leq d$ in the usual sense and its $d$th symbol
(in the usual sense) is equal to the symmetrization of $\Smbl_d(P)$ times the
identity endomorphism of $E$.

\noindent {\sl Proof:} This is true for  the case $P=D_v$ for a vector field $v$.
The general case follows from  multiplicativity. \qed

\vskip .3cm

\noindent {\bf (5.3) Classification of natural differential operators.} 
Similarly to what was done in
[Ep] [St] for the Riemannian case, it is possible to use classical invariant theory to
show, under certain assumptions, that every natural differential operator $P$
(resp. any noncommutative vector field) comes from a unique section of $\bbD_X$
(resp. $\cP_X$). These assumptions specify the way $P_{E,\nabla}$ is allowed to depend
on the connection $\nabla$ in a fixed bundle $E$. 

\vskip .2cm

Let $\Conn(E)$ be the sheaf of all $C^\infty$ connections in $E$. As well known, it is a
sheaf of torsors over the sheaf of groups $\Omega^1_X\otimes \End(E)$. For $x\in X$
let $J^p_x\Conn(E)$ be the space of all $p$-jets of connections in $E$ near $x$.
This space has a natural structure of an affine algebraic variety
over $\bbR$ (isomorphic, non-canonically, to an affine space of suitable dimension).
When $x$ varies, we get a nonlinear bundle $J^p\Conn(E)\to X$, with fibers
$J^p_x\Conn(E)$. We can speak about $C^\infty$ functions on $J^p\Conn(E)$ that are
polynomial on each fiber of this bundle, or about fiberwise polynomial morphisms of
$J^p\Conn(E)$ to another vector bundle. We denote by 
$$\Gj_p: \Conn(E)\to J^p\Conn(E) \leqno (5.3.1)$$
the ``universal differential operator'' which to each section $\nabla$ associates
the collection of its $p$-jets at each $x\in X$.

\vskip .1cm

Similarly, let $J^pE\to E$ be the bundle of $p$-jets of sections of $E$, and
$$j_p: E\to J^pE
\leqno (5.3.2)$$
be the universal differential operator of order $p$. Thus
any linear differential operator $\Phi: E\to E$ of order $p$ can be represented as
the composition, with $j_p$, of a uniquely defined linear map (morphism of
vector bundles) $\widetilde{\Phi}: J^pE\to E$. 

\vskip .2cm

\proclaim (5.3.3) Definition. A natural differential operator $P=(P_{E,\nabla})$
of order $\leq d$ is called regular, if there exists $p$ such that for each $E$
there exists a map of bundles
$$\widetilde{P} = \widetilde{P}_E: J^p\Conn(E)\times J^d E\to E,$$
polynomial in the first variable, linear in second variable
and such that for each section $s$ of $E$ and each connection
$\nabla$ in $E$  we have
$$P_{E, \nabla}(s) = \widetilde{P}(\Gj_p(\nabla), j_d(s)).$$

Informally, we require that $P_{E,\nabla}$ depends only on the $p$-jet
of $\nabla$, and in a polynomial way. \vskip .3cm

\vskip .2cm

\proclaim (5.3.4) Theorem. (a) The subalgebra in
$\cN\cD_X$ formed by regular operators, is identified with $\bbD_X$,
via the homomoprhism $h$ from (6.1.1).\hfill\break
(b) Under this identification the Lie subalgebra in $\VVect_X$ formed by
regular noncommutative vector fields, is identified with $\cP_X$.
\hfill\break
(c) The subalgebra in $\cN\cE_X$ formed by regular natural endomorphisms, is
identified with $U(\cP^\circ_X)$. \hfill\break

\noindent {\sl Proof:} We start with (c), so the $P_{E,\nabla}$ are assumed
to be endomorphisms of the bundle $E$. It is enough to work locally, so
we assume $X=\bbR^n$ with coordinates $x_1, ..., x_n$, and we also assume $E$
to be the trivial bundle of rank $N$. For any $x\in \bbR^n$ notice that
$$J^\infty_x\Conn(E):=\lim_{\longleftarrow} {} _p \,\, J^p_x\Conn(E)  = 
\Conn(E|_{\GD_{n,x}})\leqno (5.3.5)$$
is just the set of all connections in (the restriction of) $E$ on  $\GD_{n,x}$, the formal
disk near $x$. Naturality of $P$ means that the map $\widetilde{P}$
in (5.3.3) (with $d=0$) is equivariant with respect to the action of
$\Aut(E) = GL_N\, C^\infty(\bbR^n)$, the group of $C^\infty$-automorphisms
of $E$. We now apply Corollary 4.4.5 about classification of connections on
$\GD_{n,x}\simeq \GD_n$ and conclude that for any $x\in X=\bbR^n$ the operator
$$P_{E,\nabla, x}: E_x\to E_x \leqno (5.3.6)$$
given by the action of $P_{E,\nabla}$ on the fiber at $x$, depends only on
finitely many operators
$$R_{i_1, ..., i_q}(x) := \nabla_{i_1} \nabla_{i_2} ... \nabla_{i_{q-2}}
F_{i_{q-1}, i_q}: E_x\to E_x.\leqno (5.3.7)$$
Here ``depends'' means that each matrix element of $P_{E,\nabla, x}$
can be expressed as a polynomial in the matrix elements of the operators
(5.3.7).

\vskip .1cm

We now use the naturality of $P$ again, which, together with Theorem 4.4.3,
 implies that the dependence of
(the matrix elements of) $P_{E,\nabla, x}$ on (those of) the $R_{i_1, ..., i_q}(x)$
is equivariant with respect to the simultaneous action of $GL(E_x) = GL_N(\bbR)$
on all of them. 

The following statement is one of the many instances of the ``main theorem of
invariant theory'' [We].

\proclaim (5.3.8) Proposition. (a) Let $V$ be a finite-dimensional $\bbR$-vector
space of dimension $N$. Consider $\End(V)$ as an affine 
space over $\bbR$ of dimension $N^2$, and  let
$\cE_m(V)$ be the space of all polynomial maps
$$F= F(Z_1, ..., Z_m): \End(V)^m\to\End(V),$$
which are equivariant with respect to the simultaneous action of $GL(V)$. Then the
natural map
$$\bbR\langle Z_1, ..., Z_n\rangle\to \cE_m(V),$$
which takes each noncommutative polynomial into the induced map $\End(V)^m\to\End(V)$,
is surjective. 
\hfill\break
(b) In particular, under the standard embedding $\End(\bbR^N)\subset\End(\bbR^{N+1})$
any element of $\cE_m(\bbR^{N+1})$ takes $\End(\bbR^N)^m$ to $\End(\bbR^N)$, so we have
a surjection $\cE_m(\bbR^{N+1})\to\cE_m(\bbR^N)$.\hfill\break
(c) The resulting homomorphism
$$\bbR\langle Z_1, ..., Z_m\rangle\to \lim_{\longleftarrow} {} _N \,\, \cE_m(\bbR^N),$$
is an isomorphism. \qed

Applying Proposition 5.3.8 to $m$ being the number of the $R_{i_1, ..., i_q}(x)$
on which $P_{E,\nabla, x}$ depends, we conclude that for any $x\in X$
the operator $P_{E,\nabla, x}$ is represented as a noncommutative
polynomial in these $R_{i_1, ..., i_q}(x)$, and this polynomial is defined uniquely,
if we allow $N\to\infty$. This precisely means that $P$ comes from
a uniquely defined  element of the enveloping algebra of the bundle of Lie algebras $\cP^\circ_{\bbR^n}$,
since the latter is freely generated by elements that are in bijection with the $R_{i_1, ..., i_q}$. 
This proves part (c) of the theorem.

\vskip .2cm

(b). If $P$ is a regular noncommutative vector field, look at the ordinary vector field $\alpha(P)$
and its lift $D_{\alpha(P)}$ into $\VVect$. Then $P'=P-D_{\alpha(P)}$ is a regular natural endomorphism,
so it comes from a section of $U(\cP_X^\circ)$ by (c). Because $P$ and therefore $P'$ satisfies the
$\otimes$-Leibniz rule, the corresponding section of $U(\cP^\circ_X)$ is (fiberwise)
primitive, so it is really a section of $\cP^\circ_X$. 

\vskip .2cm

(a) This is reduced to (b), similarly to the argument in [St], Th. 3.7. Namely, if $P=(P_{E, \nabla})$
is a natural differential operator of order $d$, then the family of the degree $d$ symbols
(in the usual sense) of the $P_{E,\nabla}$ gives a linear system
$$\sigma_d(P_{E,\nabla}): E\to E\otimes S^dT_X$$
of natural endomorphisms parametrized by $S^dT_X$. We can work locally in coordinates $x_1, ..., x_n$,
and write
$$\sigma_d(P_{E,\nabla}) = \sum_{i_1+...+i_n=d} Q_{i_1, ..., i_n}\otimes \partial_1^{i_1} ... \partial_n^{i_n},$$
where each $Q_{i_1, ..., i_n}$ is a natural endomoprhism, which is also regular since $P$ is regular.
Lifting each $Q_{i_1, ..., i_n}$ to a section $U_{i_1, ..., i_n}$ of $U(\cP_X^\circ)$, as in (c), 
we form a section
$$\tau = \sum U_{i_1, ..., i_n} D_1^{i_1} ... D_n^{i_n}$$
of $\bbD_X$. Here $D_i = D_{\partial_i}$ is the noncommutative vector field corresponding to
$\partial_i$. We then find, using Proposition 5.2.3, that the natural differential
operator $P-h(\tau)$ has order less than $d$, and continue by induction to
argue that $P$ comes from a section $\widetilde{P}$ of $\bbD_X$. The uniqueness of $\widetilde{P}$ is
shown in a similar way: looking at the highest symbols of the $h(\widetilde{P})_{E,\nabla}$,
expanding them in the components $Q_{i_1, ..., i_n}$, and using the uniqueness already
proved for part (c). \qed

\vskip 2cm

\centerline {\bf 6. Algebro-geometric setting.}

\vskip 1cm

\noindent {\bf (6.1) Lie algebroids.} Let $k$ be a field of charcateristic 0, and 
$X$ be a scheme of
finite type over $k$ (possibly singular). Then we have a coherent sheaf $T_X$ on $X$,
called the tangent sheaf, such that $\Gamma(U, T_X) = \Der(k[U])$, if $U\subset X$
is an affine open subset. Sections of $T_X$ will be called vector fields. 
 A  Lie algebroid on $X$ is a sheaf $\cG$ of Lie-Rinehart
algebras on the Zariski topology of $X$, which is quasicoherent over $\cO_X$.
In other words $\cG$ is a quasicoherent sheaf of $\cO_X$-modules equipped with a
morphism of $\cO_X$-sheaves $\alpha: \cG\to T_X$ called the anchor map and with
a structure of a sheaf of Lie $k$-algebras such that the analog of (1.1.1-2) are
satisfied: first, $\alpha$ is a morphism of sheaves of Lie $k$-algebras, and, second,
$$[x, fy]- f[x,y] = \Lie_{\alpha(x)} (f) \cdot y. \leqno (6.1.1)$$
where $\Lie$ stands for the action of vector fields on functions (Lie derivative). 
A Lie algebroid is called transitive, if $\alpha$ is surjective.

All the constructions of (1.1-2) extend easily to this setting.
In particular, for a quasicoherent sheaf $\cM$ on $X$ we have the Atiyah algebroid
$\cA_{\cM}$, we have the concept of modules and of the enveloping algebra $U(\cG)$
of a Lie algebroid $\cG$ which is a sheaf of associative algebras on $X$. 

\vskip .2cm

\noindent {\bf (6.1.2) Example.} Assume $X$  smooth. Then  $U(T_X) = \cD_X$
is the sheaf of regular differential operators on $X$. 

\vskip .2cm

We denote $\cG^\circ=\Ker(\alpha)$. By (6.1.1) we have that
 $\cG^\circ$ is  a sheaf
of Lie $\cO_X$-algebras. 

\vskip .2cm

For any quasicoherent sheaf $\cV$ on $X$ equipped with a moprhism
$\beta: \cV\to T_X$ (anchor map), we have $\cF(\cV)$, the free Lie algebroid
generated by $\cV$, by sheafifying the construction of Theorem 2.1.2. 

Further, $\cF(\cV)= \bigcup_{d=1}^\infty \cF_{\leq d}(\cV)$ is a union of
quasicoherent subsheaves of $\cO_X$-modules (which are coherent if $\cV$
is), and this filtration is compatible with the Lie algebra structure of
sections. The sheaf of $\cO_X$-Lie algebras
$$\cF^\circ(\cV) = \Ker \{\alpha: \cF(\cV)\to T_X\} \leqno (6.3.1)$$
has the induced filtration $\{\cF^\circ_{\leq d}(\cV)\}_{d\geq 2}$
compatible with the Lie algebra structure in the fibers. Its associated
graded $\cO_X$-Lie algebra is isomorphic to $\FL_{\geq 2}(\cV)$,
the degree $\geq 2$ part of the sheaf of fiberwise free Lie algebras. 

\vskip .2cm

We denote by $\cP_X = \cF(T_X)$ the free Lie algebroid generated by $T_X$.
For each $k$-point $x\in X$ we denote
$$\cP(X,x) = \cP^\circ_X\otimes_{\cO_X} k_x\leqno (6.1.4)$$
the fiber of $\cP^\circ_X$ at $x$. This is a Lie $k$-algebra.
If $x$ is a smooth point, then $\cP(X,x)$ is free, and its
(co)homology is described as in (4.3.1-2). 

\vskip .3cm

\noindent {\bf (6.2) Noncommutative differential operators.} From now on
we assume that $X$ is a smooth algebraic variety over $k$. The constructions of
(5.1) extend easily to this setting, giving a sheaf $\bbD_X = U(\cP_X)$
of left $\cO_X$-bialgebras, which is locally free as a left or right $\cO_X$-module.
 As before, we have the two filtration on $\bbD_X$: by order and degree, with associated
graded sheaves of algebras described as in (5.2.1-2).

\vskip .2cm

As in the $C^\infty$-case, sections of $\cP_X$ and $\bbD_X$ act in bundles with
connections. There may not exist enough such bundles on all of $X$, so we have
to work with bundles defined locally. 
For every Zariski open $U\subset X$ we denote by $\Bun_\nabla(U)$ the tensor category of
vector bundles (of finite rank) on $U$ with connections. 

Let $S$ be a scheme, and $x: S\to X$ ne a morphism. Then, by (A.3-4),
$x^*\bbD_X$ (inverse image with respect to the left $\cO_X$-module structure on $\bbD_X$)
is a quasicoherent sheaf of $\cO_S$-modules. In fact, it is a sheaf of unital $\cO_S$-coalgebras.

let $P\in\Gamma(S, x^*\bbD_X)$. Then, for any open $U\subset X$, and any $(E,\nabla)\in\Bun_\nabla (U)$ we have a
morphism of sheaves on the open subscheme $x^{-1}(U)\subset S$:
$$P_{E,\nabla}: x^{-1}E\to x^*E. \leqno (6.2.1)$$
Here $x^{-1}E$ is the inverse image as a sheaf, and $x^*E$ is the inverse image as an $\cO$-module.
The construction of $P_{E,\nabla}$ is similar to the $C^\infty$-case. One starts with vector
fields $v$ on $X$ acting by operators $\nabla_v$ of covariant differentiation, and then extends
by multiplicativity (to get to $\bbD_X$) and by $\cO_S$-linearity (to get to $x^*\bbD_S = 
\cO_S\otimes_{x^{-1}\cO_X} x^{-1}\bbD_X$). 

\vskip .2cm

Let $E$ be a vector bundle on an open $U\subset X$. Then, as in (5.3.2), we have $J^pE$, the vector
bundle of $p$-jets of sections of $E$ and, for each $x: S\to X$ as above, we have a morphism
of sheaves on $x^{-1}(U)$:
$$j_{p,x}: x^{-1}E \to x^* J^pE, \leqno (6.2.2)$$
called the universal differential operator at $x$. By a $p$th order
differential operator in $E$ at $x$ we will mean 
a morphism of sheaves
$$Q: x^{-1}E \to x^*E$$
on $x^{-1}(U)$ which is a composition of a moprhism of $\cO_{x^{-1}(U)}$-sheaves
 $x^*J^pE\to x^*E$ with $j_{p,x}$.

Further, we have the sheaf $\Conn(E)$ of connections in $E$ and an affine bundle
$J^p\Conn(E)\buildrel \pi_p\over\longrightarrow U$ of $p$-jets
of connections. We consider $J^p\Conn(E)$ as an algebraic variety over $U$,
in contrast with $J^pE$ which we consider as a sheaf of $\cO_U$-modules. 
We denote by $\underline{J^p\Conn(E)}$ the sheaf of regular sections of $\pi_p$.
We have then the morphism of sheaves
$$\Gj_p: \Conn(E)\to \underline{J^p\Conn(E)}. \leqno (6.2.3)$$
Let also 
$$x^* J^p\Conn(E) := x^{-1}U \times_U J^p\Conn(E)\buildrel \pi_{p,x} \over
\longrightarrow x^{-1}(U)
\leqno (6.2.4)$$
be the induced family over $U$. 

\vskip .1cm

\proclaim (6.2.4) Definition. (a) A natural differential operator at $x$
is a system of differential operators
$$P=\{ P_{E,\nabla}: x^{-1}E \to x^* E\},$$
given for each open $U\subset X$ and each $(E,\nabla)\in\Bun_\nabla(U)$, which
is compatible with restrictions for open inclusions $U'\subset U$ and satisfies naturality for
morphisms of bundles with connections. \hfill\break
(b) A natural differential operator $P$ at $x$ is called regular of order
$\leq d$, if there exists $p\geq 0$ such that for any $U,E$ as above
there is a morphism
$$\widetilde{P} = \widetilde{P}_{U,E}: \pi^*_{p,x}J^dE \to \pi^*_{p,x} E$$
of sheaves of $\cO$-modules on $x^*J^p\Conn(E)$
with the following property. For any connection $\nabla$ in $E$ the morphism $P_{E,\nabla}$
is the composition of $j_{d,x}$ and
$$(\id\times \Gj_p(\nabla))^*(\widetilde{P}): x^*J^d E\to x^*E.$$
Here $\id\times\Gj_p(\nabla): x^{-1}(U)\to x^*J^p\Conn(E)$
is the section induced by $\Gj_p(\nabla): U\to J^p\Conn(E)$. 

\vskip .2cm

\proclaim (6.2.5) Theorem. Sections of $x^*\bbD_X$ over $S$ are in bijection with
regular differential operators at $x$. 

\noindent {\sl Proof:} Restricting to open $S'\subset S$, we see that regular operators
form a sheaf on $S_{\Zar}$, so it is enough to prove the statement
locally on $S$. For this, it is enough to work locally on $X$, so
we can assume that $X$ is an affine variety admitting an \'etale coordinate
system,
i.e., an \'etale map
$$(x_1, ..., x_n): X\to \bbA^n.$$
By Proposition 2.2.6, this allows us to trivialize $\cP_X$ and thus $\bbD_X$
as left $\cO_X$-modules, 
introducing sections $D_i\in\cP(X)$ as in Example 3.3, and writing a general element
of $\bbD(X)$ uniquely in  the form (sum with only finitely many nonzero summands):
$$\sum_{p=0}^\infty \sum_{i_1, ..., i_p=1}^n f_{i_1, ..., i_p} D_{i_1} ... D_{i_p},
\quad f_{i_1, ..., i_p}\in\cO(X).$$
Accordingly, a general section of $x^*\bbD_X$ is written uniquely in a similar form
but with $f_{i_1, ..., i_p}\in\cO(S)$. After this, the proof is achieved in the
same way as for Theorem 5.3.4, with purely notational changes.

\vskip 2cm

\centerline{\bf 7. The path groupoid.}

\vskip 1cm

\noindent {\bf (7.1) $C^\infty$ and analytic cases.}
Here is a convenient way of defining an equivalence relation on the
space of parametrized paths so as to get a groupoid.

\proclaim (7.1.1) Definition. (cf. [HKK])
 Let $X$ be a $C^\infty$-manifold, and
$\gamma, \gamma': [0,1]\to X$ be two piecewise smooth parametrized paths with the
same beginning $x=\gamma(0)=\gamma'(0)$ and end $y=\gamma(1) = \gamma'(1)$.
A piecewise smooth homotopy $\sigma: [0,1]^2\to X$ between $\gamma$ and $\gamma'$ will be
called  thin, if, for any $(\lambda,\mu)\in [0,1]^2$ where
$\sigma$ is smooth, the differential $d_{(\lambda, \mu)}\sigma$ has rank $\leq 1$.

Let $\Pi_X$ be the set of piecewise smooth paths in $X$ modulo thin homotopies, and
$\Pi_X(x,y)$ be the set of classes of paths with $\gamma(0)=x$ and $\gamma(1)=y$.
Then $\Pi_X$ is (the set of morphisms of) a groupoid with $X$ as the set of objects:
thin homotopies account for both reparametrizations and cancellations of any segment
followed by the same segment run in the opposite direction. 
We call $\Pi_X$ the {\em path groupoid} of $X$. 

\vskip .2cm

It is natural to view $\cP_X$ as the Lie algebroid of $\Pi_X$. Indeed, for $x\in X$ 
consider the tensor functor
$$\Phi_x: (\Bun_\nabla(X), \otimes) \to (\Vect_{\bbR},\otimes) , \quad (E,\nabla)\mapsto E_x,\leqno (7.1.2)$$
associating to a vector bundle its fiber at $x$.
 Then, for any $\gamma\in\Pi_X(x,y)$ and any $(E,\nabla)\in\Bun_\nabla(X)$,
we have the holonomy operator $H_{E,\nabla}(\gamma): E_x\to E_y$, and these
operators give a natural transformation of functors $\Phi_x\to\Phi_y$ compatible
with the tensor functor structures. On the other hand, by (3.2), sections of $\cP_X$
evaluated at any $x$, give similar transformations except they satisfy 
the $\otimes$-Leibniz rule which is the infinitesimal version of compatibility with
the tensor product. 

Thus, sections of $\cP_X$ can be seen as vertical vector fields on
$$\Pi_X\buildrel s\over\longrightarrow X,\leqno (7.1.3)$$
where $s$ is the source map (beginning of the path), see (A.1.3). 

\vskip .2cm

 Let now  $X$ be a complex analytic manifold.
A homotopy $\sigma$ as above will be called {\em $\bbC$-thin}, if for any $(\lambda,\mu)$
where $\sigma$ is smooth, the image of $d_{(\lambda, \mu)}\sigma$ is contained in a
1-dimensional $\bbC$-linear subspace of $T_{\sigma(\lambda,\mu)}X$. We denote by
$\Pi_X^\bbC$ the set of piecewise smooth paths in $X$ modulo $\bbC$-thin homotopies.
Thus $\Pi_X^\bbC$ is a groupoid, a quotient of $\Pi_X$. 

Given a holomorphic vector bundle $E$ with a holomorphic connection $\nabla$,
the holonomy operators along the paths are unchanged under $\bbC$-thin homotopies, and
so  give rise to a functor
$$H_{E,\nabla}: \Pi_X^\bbC\to\Vect_\bbC.\leqno (7.1.4)$$

\vskip .3cm

\noindent {\bf (7.2) Kontsevich spaces.} Let $k$ be a field of
characteristic 0, and  $X$ is a smooth projective
 algebraic variety over $k$. As in [Ma], \S V.1.4, let
$B(X) \subset \Hom({\operatorname{Pic}}(X), \bbZ)$ be the set of
homomorphisms which are $\geq 0$ on each ample line bundle $L$. 
 Let $\beta\in B(X)$, and $M_\beta = \overline{M}_{0,2}(X,\beta)$
be the moduli stack of stable 2-pointed rational curves in $X$ of degree $\beta$.
Recall ([Ma], \S V.3.2) that for any $\bbC$-scheme $S$ the $S$-points of $M_\beta$
form a groupoid $M_\beta(S)$ whose objects are data $(C,x,y,f)$ where:

\vskip .2cm

\noindent (7.2.1) $C\buildrel\pi\over\longrightarrow S$ is a flat family of proper
curves with every geometric fiber being a union of $P^1$'s whose intersection
graph is a tree. 

\vskip .1cm

\noindent (7.2.2) $x,y: S\to S$ are two sections of $\pi$ which are everywhere
disjoint and whose values at any geometric point $s\in S$ 
are smooth on $\pi^{-1}(s)$. 

\vskip .1cm

\noindent (7.2.3) $f: C\to X$ is a morphism of schemes such that for any $s\in S$
as above the degree of $f|_{\pi^{-1}(s)}$ is $\beta$. 

\vskip .1cm

These data are required to satisfy the stability condition ([Ma] \S V.1.2).
Morphisms between $(C,x,y,f)$ and $(C',x', y', f')$ are isomorphisms $C\to C'$
preserving all the data. 

\vskip .2cm

The disjoint union  $M=\coprod M_\beta$ can be seen as an algebro-geometric
analog of the path groupoid. Indeed, we have the projections
$$s,t: M_\beta\to X, \quad (C,x,y,f)\mapsto x,y.\leqno (7.2.4)$$
Further, gluing  the mapped curves together followed  by the stabilization morphism
discussed in [Ma], \S V.1.7, gives morphisms of  stacks 
$$m_{\beta,\beta'}: M_\beta \times_X M_{\beta'}\to M_{\beta+\beta'} \leqno (7.2.5)$$
remindful of the composition map (2.1.2) of a groupoid. However, the $m_{\beta, \beta'}$
do not possess inverses or unit. Still, the following statement 
expresses the similarity between the two objects.

\proclaim (7.2.6) Proposition. Let $k=\bbC$. Then there
 exist  natural morphisms $u_\beta: M_\beta(\bbC)\to \Pi_X^\bbC$
which take the compositions (7.2.5) into the composition of the groupoid $\Pi_X^\bbC$.

\noindent {\sl Proof:} Let $(C,x,y,f)\in M_\beta(\bbC)$. Thus $C$ is a curve over $\bbC$
which is a union of a tree of projective lines. In particular, $C$ is simply connected.
Next, $x,y$ are  two smooth points of $C$ and $f: C\to X$ is a map of degree $\beta$. 
Since $C$ is simply connected,  $x$ and $y$ can be joined by
a path in $C$ which is unique up to homotopy. After applying $f$ any such homotopy
will be a $\bbC$-thin homotopy in $X$. Thus we have a unique  $\bbC$-thin homotopy
class of paths in $X$ joining $f(x)$ and $f(y)$. This defines $u_\beta$. The rest is clear. \qed

\vskip .2cm

\noindent {\bf (7.2.7) Remark.} Similarly, $\bbC$-points of
 the stack $\overline{M}_{0,n}(X,\beta)$
of stable $n$-pointed rational curves, $n>2$, can be mapped into $B_{n-1}\Pi_X^\bbC$,\
the $(n-1)$-st component of the simplicial classifying space of the groupoid $\Pi^\bbC_X$. 

\vskip .3cm

\noindent {\bf (7.3) The formal path groupoid.}
 Let $k$ be a field of characteristic $0$, and
$X$ be a smooth algebraic variety over $k$ (not necessarily projective).
 Denote by  $\widehat{\Pi}_X = e^{\cP_X}$ be the formal groupoid integrating the Lie algebroid $\cP_X$,
as constructed in (A.6). Thus its scheme of objects is $X$, and the ind-scheme of morphisms
(for which we retain the notation $\widehat{\Pi}_X$) is
$$\widehat{\Pi}_X = \Spf(\GA) = ``\lim_{\longrightarrow} {''} _d \,\, \Spec(\GA_d), \quad
\GA_d = \shHom_{\cO_X}(\bbD_X^{\leq d}, \cO_X). \leqno (7.3.1)$$
We call $\widehat{\Pi}_X$ the {\em formal path groupoid} of  $X$.

\vskip .2cm

\noindent {\bf (7.3.2) Remark.}
 For $k=\bbC$ the ind-scheme $\widehat{\Pi}_X$ 
 provides
an algebro-geometric
model for the formal neighborhood of $X$ in $\Pi_X^\bbC$. This can be understood
as follows. By Theorem A.6.5, a morphism from a scheme $S$ into $\widehat{\Pi}_X$
is the same as a morphism $y: S\to X$ together with a groupoid-like section $g$ of $y^*\bbD_X$.
On the other hand, suppose that $y$ is a $\bbC$-point of $X$, and $\gamma$ be a path in $X$
with endpoint $y$. Assuming that $\gamma$ lies in a coordinate patch on $X$
with coordinates $x_1, ..., x_n$ (vanishing at $y$), we can associate to
$\gamma$ the formal series
$$E_\gamma = \sum_{p=0}^\infty \sum_{i_1, ..., i_p=1}^n \biggl(\int_\gamma^{\rightarrow} dx_{i_1}
... dx_{i_p}\biggr) D_{i_1} ... D_{i_p}, \leqno (7.3.3)$$
which can be seen as lying in the formal completion of the fiber of $\bbD_X$ at $x$.
Here the $\int_\gamma^{\rightarrow}$ are the iterated integrals of Chen [C], and his fundamental
shuffle relations imply that $E_\gamma$ is, formally, a group(oid)-like element.
Indeed, the $D_i$, being sections of $\cP_X$, are primitive. Although the concept of
constant differential forms $dx_i$ and thus of their iterated integrals depends on the choice
of coordinates, the transformation rules for sections of sheaves $\cP_X$ and $\bbD_X$ are
such that the generating series $E_\gamma$ is invariantly defined.

\vskip .2cm

We now extend this point of view to relate $\widehat{\Pi}_X$ to the Kontsevich
moduli spaces by providing an algebraic version of Proposition 7.2.6. 
We need to introduce some terminology.

Let $S$ be a reduced scheme, and $(C,x,y,f)\in M_\beta(S)$. Then $C$ is reduced as well.
Denote by $C^0\subset C$ the minimal closed fiberwise connected relative subcurve
containing $x(S)$ and $y(S)$. 
We say that $(C,x,y,f)$ is {\em contracting}, if $f$ is constant on the fibers of the
projection $\pi^0: C^0\to S$, i.e., $f|_{C^0}$ factors through a morphism $S\to X$. 

Let now $S$ be arbitrary, and $(C,x,y,f)\in M_\beta(S)$. We then have $C_{\red}=\pi^{-1}(S_{\red})$,
so we have a stable  $\pi_{red}: C_{\red}\to S_{\red}$ with sections $x_{\red}, y_{\red}$
and map $f_{\red}: C_{\red}\to X$ given by restriction. They form an $S_{\red}$-point of $M_\beta$.
We say that $(C,x,y,f)$ is {\em almost contracting}, if $(C_{\red}, x_{\red}, y_{\red}, f_{\red})$
is contracting. Let $\widehat{M}_\beta\subset M_\beta$ be the substack whose $S$-points
are almost contracting $(C,x,y,f)$. We can think of $\widehat{M}_\beta$ as a certain formal neighborhood
in $M_\beta$. 

Note that the compositions (7.2.5) restrict to
$$\widehat{m}_{\beta,\beta'}: \widehat{M}_\beta\times_X
 \widehat{M}_{\beta'}\to\widehat{M}_{\beta+\beta'}.\leqno (7.3.4)$$

\proclaim (7.3.5) Theorem. There exist morphisms of stacks $\widehat{u}_\beta: \widehat{M}_\beta\to
\widehat{\Pi}_X$ which are compatible with the source, target and the composition maps.

\noindent {\sl Proof:} To construct $\widehat{u}_\beta$ means, by definition, to construct,
for any scheme $S$, a map $\widehat{u}_{\beta, S}: \widehat{M}_\beta(S)\to\widehat{\Pi}_X(S)$
 constant on isomorphism classes in $\widehat{M}_\beta(S)$, in a way compatible with the
base change for the {\em fppf} topology (on which the stacks are defined).

Since $M_\beta$ (and thus $\widehat{M}_\beta$) is a stack, and $\widehat{\Pi}_X$,
being an ind-scheme, is a sheaf of sets on the {\em fppf} topology, we can assume $S$ affine.
Since $M_\beta$ is a stack of finite type over $k$, we can assume $S$ be be of finite type as well. 

Let $p=(C,x,y,f)$ be an $S$-point of $\widehat{M}_\beta$. The fact that it is almost
contractible implies that the $S$-points $fx, fy: S\to X$ are infinitesimally close, i.e.,
coincide on $S_{\red}$. Let $U\subset X$ be Zariski open, and $(E,\nabla)\in\Bun_\nabla(U)$. 
Because $fx, fy$ are infinitesimally close, we have the morphism of sheaves (restriction)
$$\Res: (fy)^{-1} E\to (fx)^*E. \leqno (7.3.6)$$
On the other hand, by Theorem A.6.5, an element of $\widehat{\Pi}_X(S)$ is the same as a pair
$(z,g)$, where $z: S\to X$, and $g\in\Gamma(S, z^*\bbD_X)$ is a groupoid-like element.
By Theorem 6.2.5, a section $g$ of $z^*\bbD_X$ over $S$ is the same as a system of
morphisms of sheaves
$$g_{E,\nabla}: z^{-1}E \to z^*E, \leqno (7.3.7)$$
given for any $U, E, \nabla$ as above and forming a regular natural differential
operator along $z$. For $g$ to be  groupoid-like is equivalent to the fact that the
$g_{E,\nabla}$ commute with tensor products of bundles with connections. 

We are going to construct there two types of data as follows. Given $p\in\widehat{M}_\beta(S)$
as above, we set $z=fy$. We then construct $g_{E,\nabla}$ as the composition of $\Res$
with a morphism (``holonomy'')
$$H_{E,\nabla}(p): (fx)^*E \to (fy)^*E.\leqno (7.3.8)$$
In other words, we prove the following fact which will imply the existence of $\widehat{u}_\beta$. 

\proclaim (7.3.9) Theorem. (a) For each $(E,\nabla)\in\Bun_\nabla(U)$ as above there exists
an isomorphism of $\cO_S$-sheaves $H_{E,\nabla}(p)$ as in (7.3.8), compatible
with restrictions to $\widetilde{U}\subset U$ and fppf base change for $\widetilde{S}\to S$.
\hfill\break
(b) The composition $H_{E,\nabla}(p)\circ \Res$ is a regular natural differential operator.
\hfill\break
(c) For fixed $p$ the $H_{E,\nabla}(p)$ are compatible with the tensor product.

Note that it enough to assume (by restricting $U$ if necessary) that the embedding $U\hookrightarrow X$
is an affine morphism. 

\vskip .3cm

\noindent {\bf (7.4) Proof of Theorem 7.3.9.} Consider the open subschemes
$$C_U = f^{-1}(U)\subset C, \quad S_U= (fx)^{-1}(U)=(fy)^{-1}(U)\subset S.$$
By restricting $S$ if necessary, we can assume that $S_U=S$. Let $\pi_U: C_U\to S$
be the restriction of $\pi$, and $f_U: C_U\to X$ be the restriction of $f$. We have then the
bundle $f_U^*E$ on $C_U$ and a relative connection
$$\nabla_{C/S}: f_U^*E\to \Omega^1_{C/S}\otimes f_U^*E. \leqno (7.4.1)$$
Note that $C$ is of finite type because we assumed $S$ to be of finite type.
Let $\widehat{C}$ be the formal neighborhood of $C^0_{\red}$ in $C$ (or, what is the same,
in $f^{-1}(U)$). We consider it as a topologically ringed space $(C^0_{\red})_{\Zar}$
with sheaf $\cO_{\widehat{C}}$. We denote by $\widehat{\pi}:\widehat{C}\to S$
the projection and for each quasicoherent sheaf $\cF$ on $C$ write
$$\Gamma(\widehat{C},\cF) = \Gamma(\cF\otimes\cO_{\widehat{C}}).$$
For any geometric point $s\in S$ the preimage $\widehat{\pi}^{-1}(s)$
is a formal scheme over the field $K=k(s)$ consisting of
 $\pi^{-1}(s)\cap C_0$  (a chain of $P^1$'s),  to which are attached several
``tails'' at the points of intersection with other components of $\pi^{-1}(s)$.
Each such tail is isomorphic to $\Spf \,\, K[[t]]$.

Theorem 7.3.9 is a consequence of the following fact.

\proclaim (7.4.2) Lemma. There exists a unique
$$\Phi_x\in\Gamma(\widehat{C}, \shHom(\pi_U^*x^*f_U^*E, f_U^*E))$$
(``fundamental solution''), satisfying 
$$\nabla_{C/S}\Phi_x =0, \quad x^*\Phi_x = \id\in\Hom_S(x^*f^*_UE, x^*f^*_UE). $$

Indeed, we then define
$$H_{E,\nabla}(p) = y^*\Phi_x\in\Hom_S(x^*f_U^*E, y^*f_U^*E).\leqno (7.4.3)$$
Its multiplicativity in tensor products as well as naturality in restrictions and
base changes follow from the uniqueness, so parts (a) and (c) of the theorem
would be proved. The validity of part (b), i.e., the fact that
$$H_{E,\nabla}(p)\circ\Res: y^{-1}f^{-1}E \to y^* f^*E$$
depends on sections of $y^{-1}f^{-1}E$ and on $\nabla$ only through their finite jets,
will follow from the construction of $\Phi_x$ which we are about to give.

\vskip .3cm

\noindent {\bf (7.4.4) Construction of $\Phi_x$.}  We proceed by
``induction on the order of nilpotency'' of $S$ (assumed to be of finite type).
First, consider the case when $S$ is reduced, so $S$ and $C$ are algebraic
varieties over $k$. Then $C^0\subset C$ is a union of irreducible components,
and the restriction of $f$ to $C^0$ factors through $S$. In particular, $f_U^*E$
is trivial over the fibers of $C^0\to S$ and the connection $\nabla_{C/S}$
is also trivial over such fibers. In this case both the existence and the
uniqueness of $\Phi_x$ is clear. Indeed, over $C^0$ it is (and must be)  equal to
the identity. Further, on the formal neighborhood of $C^0$ in $C$ (which is $\widehat{C}$)
it is uniquely extended by solving, over each geometric point $s\in S$, finitely many
 initial value problems for 
 differential equations in the ring of formal power series $k(s)[[t]]$.
These  equations  correspond
to the  ``tails'' attached, as above, at the 
 intersection point of $(\pi^0)^{-1}(s)$ with
irreducible components of $\pi^{-1}(s)$ not in $C^0$.  

\vskip .2cm

For the inductive step, we consider an arbitrary $S$ of finite type and a closed subscheme
$S'\subset S$ whose sheaf of ideals $I_{S'}\subset \cO_S$ satisfies $I_{S'}\cdot\sqrt{\cO_S} =0$. 
This means that $I_{S'}$, as a coherent sheaf on $S$, is supported, scheme-theoretically, on $S_{\red}$.
Let $C'=\pi^{-1}(S')$. Then the restriction $\pi': C'\to S'$ together with $f' = f|_{C'}$,
as well as $x'=x|_{S'}, y'=y|_{S'}$ form an $S'$-point $p'$ of $\widehat{M}_\beta$.
In particular,  we have the formal
scheme $\widehat{C}'$, the formal neighborhood of $C^0_{\red}$ in $C'$ etc. 
  
\vskip .2cm

Let $I_{C'}\subset \cO_C$ be the sheaf of ideals of $C'$. Because $C/S$ is flat, we have
$$I_{C'} = \pi^* I_{S'}, \quad I_{C'}\cdot\sqrt{\cO_C}=0.\leqno (7.4.5)$$
Further, flatness is retained by the completion $\widehat{C}\to S$,
so the sheaf of ideals $I_{\widehat{C}'}\subset \cO_{\widehat{C}}$ satisfies the
analog of (7.4.5). In addition, since $C$ is a curve of relative arithmetic genus 0,
for any coherent sheaf $\cF$ on $S$ we have
$$R^1\widehat{\pi}_* (\widehat{\pi}^* \cF)=0.\leqno (7.4.6)$$

\vskip .2cm

We assume by induction  Lemma 7.4.2 to be proved for $\widehat{C}'$. So let
$$\Phi'=\Phi_{x'}\in\Gamma(\widehat{C}'. \Hom(\pi_U^*(x')^*(f'_U)^*E, (f'_U)^*E))
\leqno (7.4.7)$$
satisfies the conditions of the lemma. Denote by
$$\cH \subset \shHom_{\widehat{C}'}(\pi_U^* x^*f_U^*E, f_U^*E)
\leqno (7.4.8)$$
the sheaf of all $\Phi$ restricting to $\Phi'$ over $\widehat{C}'$.
This is a torsor over the sheaf of all $\Phi$ restricting to 0, i.e., over
$$\cA= \shHom_{\widehat{C}}(\pi_Ux^*f_U^*E, f_U^*E\otimes I_{C'}).\leqno (7.4.9)$$
We can assume, by restricting $U$ if necessary, that $E$ is trivial as a vector bundle
on $U$. Then $\cA$, which is, by the properties of  $I_{C'}$,
 a coherent sheaf on $C_{\red}$, is trivial on fibers of
$C_{\red}\to S_{\red}$, so its first direct image is 0 by (7.4.6).
This means that locally on $S$ the torsor $\cH$ has a global section.
Since we work locally on $S$ anyway, we can assume that $\cH$ is trivial as a torsor. 
So there exists $\Phi_0$ extending $\Phi'$. We now modify $\Phi'$ be a section
$\Psi$ of $\cA$ so as to satisfy the conditions of Lemma 7.4.2. 
Without such modification we have only that
 $$\omega := \nabla_{C/S}(\Phi')\in\Gamma (\widehat{C}, 
\shHom(\pi_U^*x^*f_U^*E, f_U^*E)\otimes
 \Omega^1_{C/S}\otimes I_{C'} ).\leqno (7.4.10)$$
As $I_{C'}$, and thus $\omega$, is supported on $C_{\red}$, we see that the restriction of $\omega$  
to each fiber of $C_{\red}'\to S_{\red}$  must vanish since each  such  fiber is a union of $P^1$'s.
So modification of $\Phi'$ is not necessary over $C_{\red}'$, and we need to look for $\Psi$ as
a section of $\cA$  which
vanishes on $C_{\red}'$. As before, this is done by 
 solving  a differential equation in the ring of formal power series
over each geometric point $s\in S$ and each intersection point of $(\pi^0)^{-1}(s)$ with other 
irreducible components.

The uniqueness is proved by the same inductive reasoning, using, at each step, the fact
that the solution of an initial value problem for
a differential equation in the ring of formal power series is unique.
This finishes the proof of Lemma 7.4.2 and thus of Theorems 7.3.9 and 7.3.5.

\vskip 2cm

\centerline {\bf 8. Remarks and further directions.}

\vskip 1cm

\noindent {\bf (8.1) Integral kernels on $\Pi_X$.} 
Of great interest are various versions of the convolution algebra of $\Pi_X$,
see [Co], \S II.5 for general background. 
Naively, this algebra should consist
of fiberwise measures on (7.1.3), i.e., of ``kernels'' $K=K(x,\gamma)\cD\gamma$
defined for $x\in X$, $\gamma\in s^{-1}(x)$, which behave as functions in $x$
and measures in $\gamma$. At this naive level, such a $K$ gives, for each $(E,\nabla)\in\Bun_\nabla(X)$,
an operator $P^K_{E,\nabla}$ in sections $f$ of $E$ (Feynman-Kac formula):
$$(P^K_{E,\nabla} f)(x) = \int_{\gamma\in s^{-1}(x)} K(x,\gamma)H_{E,\nabla}(\gamma)\bigl(f(t(\gamma))\bigr)\cD\gamma.
\leqno (8.1.1)$$
Here $t(\gamma)$ is the endpoint of $\gamma$ and $H_{E,\nabla}(\gamma): E_x\to E_{t(\gamma)}$ is
the holonomy of $\nabla$ along $\gamma$. 
. For example, ``distribution kernels'' obtained as
iterated derivatives of the delta function on $X\subset\Pi_X$, give rise to $P^K_{E,\nabla}$ being
the natural differential operators corresponding to  sections of $\bbD(X)$, as in (6.1).

Formula (8.1.1) involves path integration and so may be difficult to justify rigorously.
On the other hand, the result, considered for all $(E,\nabla)$, is manifestly natural
with respect to morphisms in $\Bun_\nabla(X)$. Extending the approach of (5.1), one
can look at natural (systems of) operators $\{P_{E,\nabla}: E\to E\}$ given in all bundles
with connections on $X$, but of more general nature: pseudo-differential, Fourier integral
operators etc. By the above, various algebras of such natural operators provide
regularized versions of the convolution algebra of $\Pi_X$. Intuitively, 
knowledge of $P^K_{E,\nabla}$ for all $(E,\nabla)$ determines $K$ uniquely. 

\vskip .2cm

\noindent {\bf (8.1.2) Example.} Let $X$ be a Riemannian manifold, and $\Delta = \Delta_X
\in\bbD(X)$ be its noncommutative Laplacian (5.1.4). We have then a natural operator $\Theta$
with $\Theta_{E,\nabla} = \exp(-\Delta_{E,\nabla})$ being the heat operator. By the result of
Bismut [Bi], $\Theta$ corresponds, as in (7.1.3), to the Wiener (Brownian motion) measure
on $\Pi_X$. To be precise, this interpretation requires completing $\Pi_X$ to include continuous
paths and understanding $H_{E,\nabla}(\gamma)$ by using stochastic differential equations.
Cf. [Kap2] for the discussion of the flat case $X=\bbR^n$. 

\vskip .3cm

\noindent {\bf (8.2) Algebraic theory of $\bbD$-modules.} Let $X$ be a smooth
algebraic variety. By analogy with sheaves of modules over $\cD_X$, the
sheaf of usual differential operators [Kas], one can study sheaves of modules
over $\bbD_X$, quasicoherent over $\cO_X$, in a purely algebraic fashion.
Such modules provide a language for studying connections with singularities.
On the other hand, sections of $\bbD_X$ are interpreted as translation
invariant differential operators on $\Pi_X$, see (7.1.3), so $\bbD_X$-modules
describe some systems of linear PDE on the space of paths. 

\vskip .3cm

\noindent {\bf (8.3) Quotients of $\cP_X$ and geometric structures.} 
Many differential-geometric structures on $X$ can be formulated in terms of
appropriate quotients of $\cP_X$. For example, a Riemannian metric on $X$
specifies the subcategory of bundles with connections satisfying the Yang-Mills equations
and thus gives rise to a quotient groupoid $\cP_X^{\operatorname{YM}}$ of $\cP_X$.
If $\dim(X)=4$, we have a further quotient $\cP_X^+$ governing self-dual connections.
These are curvilinear versions of the (self-dual) Yang-Mills algebras of Nekrasov
and Connes-Dubois-Violette [Ne] [CDV]. They will be studied in a subsequent paper. 

\vskip 2cm

\centerline {\bf Appendix. Formal integration of Lie algebroids.}

\vskip 1cm

\noindent {\bf (A.1) Lie groupoids.} Recall that a { groupoid} $\GG$ is a category in which
all morphisms are isomorphisms. 
 We will write $\GG_0= \Ob(\GG)$ for the class of objects
and $\GG_1 = \Mor(\GG)$ for the class of morphisms of $\GG$. In the sequel we assume that
$\GG$ is small, i.e., that the $\GG_i$  are sets. Thus we have maps
$$s,t: \GG_1\to \GG_0, \quad e: \GG_0\to \GG_1, \leqno (A.1.1)$$
where $s,t$ associate to a morphism its source and target objects while
$e$ associates to an object its identity (unit) moprhism. In addition, the composition
of morphisms can be seen as a map
$$m: \GG_1\times_{\GG_0}\GG_1\to \GG_1,\leqno (A.1.2)$$
where the fiber product is taken with respect to the pair of maps $s,t: \GG_1\to\GG_0$. 

\vskip .2cm

One can speak about groupoid object in any category $\cC$ with fiber products.
Taking $\cC$ to be the category of $C^\infty$-manifolds, we get a concept of
a Lie groupoid. Thus a 
{Lie groupoid} is a groupoid $\GG$ in which both $\GG_0$ and $\GG_1$ are equipped with
structures of $C^\infty$-manifolds, and the maps (A.1.1-2) are smooth. When $\GG_0$
is a point, a Lie groupoid is the same as a Lie group. 
 See [Mack] for more
background and examples of  Lie groupoids. 

\vskip .2cm

We recall the fundamental construction which to any Lie groupoid $\GG$ associates
a Lie algebroid $\cG = \Lie(\GG)$. Denote the manifold of objects by $X=\GG_0$. 
Then, as a vector bundle,
$$\cG  =  e^* (T_{\GG_1/ \GG_0}), \leqno (A.1.3)$$
where the relative tangent bundle of $\GG_1/\GG_0$ is taken with respect to the
projection $s$. In fact, sections of  $\cG$ can be seen as vertical vector
fields on $\GG_1\buildrel s\over\longrightarrow \GG_0$ which are right
invariant with respect to the composition map, and the Lie bracket
of such vector fields makes $\cG$ into a sheaf
of Lie algebras. Further, the differential of $t: \GG_1\to\GG_0$
defines the anchor map $\alpha: \cG\to T_X$, and these structures make
$\cG$ into a Lie algebroid. 
 See [Mack] for more details. 

\vskip .2cm

It is known that any finite-dimensional Lie algebra can be integrated to a Lie
group. The corresponding problem for Lie algebroids is higly nontrivial, 
see [Mack], and the integration is not always possible. 
In the following  we present
a groupoid version of an easier result: that any Lie algebra can be integrated
to a {\em formal} (Lie) group. 

\vskip .3cm

\noindent {\bf (A.2) (Formal)  ind-schemes.}  We follow the same conventions on ind-schemes as in [KV1-2],
see also [De] for general background on ind-objects. 

\vskip .1cm

Thus, by an ind-scheme $\cX$ over $k$ we mean a formal inductive limit
$\cX = ``\lim\limits_{\longrightarrow}{''} {}_{i\in I} X_i$, where $I$ is a filtered
poset, and $(X_i)_{i\in I}$ is an inductive system of $k$-schemes indexed by $I$.
It is further assumed that the structure maps $X_i\to X_j$, $i\leq j$, are closed
embeddings.  An ind-scheme $\cX$ can be identified with the corresponding (ind-)representable
functor $h_\cX$ on the category of $k$-schemes:
$$h_\cX(S) = \lim_{\longrightarrow} {} _{i\in I}\,\, \Hom(S, X_i). \leqno (A.2.1)$$
This means that if $\cY = ``\lim\limits_{\longrightarrow}{''}{} _{j\in J} Y_j$ is
another ind-scheme (with a possibly different indexing set $J$), then morphsisms
$\cX\to \cY$ are the same as natural transformations of functors $h_\cX\to h_{\cY}$. 

\vskip .2cm

\noindent {\bf (A.2.2) Examples.} (a) Let $A$ be a commutative  topological $k$-algebra
which can be represented as a projective limit $A=\lim\limits_{\longleftarrow} {} _{i\in I}\, A_i$
of discrete $k$-algebras. Then we have the ind-scheme
$$\Spf(A) = ``\lim_{\longrightarrow}{''} _{i\in I} \Spec(A_i).$$
In particular, we will use the {\em $n$-dimensional formal disk} over $k$ which is
the ind-scheme
$$\GD_n = \Spf\,\, k[[x_1, ..., x_n]] = ``{\lim_{\longrightarrow}{}}{ ''} {} _d \,\, \Spec\,\, k[x_1, ..., x_n]/
(x_1, ..., x_n)^d.$$

\vskip .1cm

(b) Let $Y$ be a scheme, and $X\subset Y$ be a closed subscheme with sheaf of ideals $I_X\subset \cO_Y$.
Then we have the ind-scheme
$$X^{(\infty)}_Y = ``\lim_{\longrightarrow}{''} _d \,\, X^{(d)}_Y, \quad X^{(d)}_Y = 
 \Spec \, (\cO_Y/I_X^{d+1}),$$
called the formal neighborhood of $X$ in $Y$. The scheme $X^{(d)}_Y$ is called the $d$th infinitesimal
neighborhood of $X$ in $Y$. 

\vskip .2cm

For a scheme $X$ we denote by $X_{\red}\subset X$ the maximal reduced (nilpotent-free) subscheme. 
We then extend this to ind-schemes, putting $\cX_{\red} = ``\lim\limits_{\longrightarrow}{''} X_{i, \red}$,
if $\cX= ``\lim\limits_{\longrightarrow}{''} X_i$. 

\vskip .1cm

\proclaim (A.2.3) Definition. An ind-scheme $\cX$ will be called formal, if $\cX_{\red}$ is
an ordinary scheme. 

Thus a formal ind-scheme can be seen as a limit of nilpotent extensions of a scheme. 
This concept is more flexible than that of formal schemes as defined in [Gr].

\vskip .2cm

For a scheme $X$ we denote by $X_{\Zar}$ the underlying topological space of $X$ with the
Zariski topology, so $X$ is a ringed space $(X_{\Zar}, \cO_X)$. Following Haboush [Hab], we associate to
an ind-scheme $\cX= ``\lim\limits_{\longrightarrow} {''}_{i\in I}  X_i$ a topological space $\cX_{\Zar}$ with a sheaf of
topological rings $\cO_{\cX}$. Explicitly,
$$\cX_{\Zar} = \lim_{\longrightarrow} {} _{i\in I}\,\,  X_{i, \Zar}, \leqno (A.2.4)$$
(inductive limit in the category of topological spaces), while
$$\cO_{\cX} = \lim_{\longleftarrow} {} _{i\in I}\,\,  \epsilon_{i, *} \cO_{X_i}, \leqno (A.2.5)$$
where $\epsilon_i: X_{i, \Zar}\to X_{\Zar}$ is the natural inclusion. It was proved in 
{\em loc. cit.} that the correspondence $\cX\mapsto (\cX_{\Zar}, \cO_\cX)$ embeds the
category of ind-schemes
as a full subcategory of the category of locally topologically ringed spaces.

\vskip .2cm

\noindent {\bf (A.2.6) Example.} In the situation of Example A.2.2(b) we have
$$(X^{(\infty)}_Y)_{\Zar} = X,  \quad \cO_{X^{(\infty)}_Y} = \widehat{\cO}_{Y, X} := \lim_{\longleftarrow} \,
\cO_Y/I_X^d.$$

\vskip .2cm

\proclaim (A.2.7) Definition. 
By a { formal groupoid}
 we will mean  a groupoid 
object $\GG = (\GG_0, \GG_1)$ in the category of formal ind-schemes over $k$
such that each of the morphisms  $e, s, t$ induces an isomorphism
between  $(\GG_0)_{\red}$ and $(\GG_1)_{\red}$. A formal group is a formal
groupoid with $\GG_0=\Spec(k)$. 

\vskip .2cm

\noindent {\bf (A.2.8) Example.} Let $\Gg$ be a Lie $k$-algebra and $U=U_k(\Gg)$ be its
universal enveloping algebra. As well known, $U$ is a cocommutative Hopf algebra with
comultiplication $\Delta$ defined uniquely by setting each $x\in\Gg\subset U$ to be
primitive:
$$\Delta(x) = x\otimes 1 + 1\otimes x. \leqno (A.2.9)$$
Let $U^\vee = \Hom_k(U,k)$ be the linear dual of $U$ as a vector space. 
Then $\Delta$ makes $U^\vee$ into a topological commutative algebra, and
$$e^\Gg:= \Spf(U^\vee)$$
is a formal group. If $\Lambda$ is any commutative $k$-algebra, then $\Lambda\otimes_k U$
is a Hopf $\Lambda$-algebra, and one has
$$\Hom(\Spec(\Lambda), e^\Gg) = (\Lambda\otimes_k U)^{\gr}.\leqno (A.2.10)$$
Here the subscript ``$\gr$'' means the set of invertible group-like elements, i.e.,
invertible elements $g$ satisfying $\Delta(g) = g\otimes g$. It is this approach
that we generalize to Lie algebroids. 

\vskip .3cm

\noindent {\bf (A.3) Left bialgebras.} Let $A$ be a commutative $k$-algebra,
included into a possibly noncommutative algebra $H$. We do not assume $A$ to be
central in $H$. Thus there are two commuting $A$-module structures on $H$, which
we call left and right. Since $A$ is commutative, each of them can be used to form
tensor products over $A$. For examples, there are four possibilities for $H\otimes_A H$.

\vskip .1cm

In the sequel we will always understand $H\otimes_A H$ to be taken with respect
to the left $A$-module structures on both factors, i.e., to be the quotient of
$H\otimes_k H$ by the relations
$$au_1\otimes u_2 = u_1\otimes au_2, \quad a\in A, u_1, u_2\in H.
\leqno (A.3.1)$$
With this definition, $H\otimes_A H$ has one, ``left'', $A$-module structure, with
$a(u_1\otimes u_2)$ given by (A.3.1) and two ``right'' structures, each commuting
with the left one:
$$(u_1\otimes u_2)a = u_1a\otimes u_2, \quad {\operatorname{resp.}}\quad
u_1\otimes u_2 a.\leqno (A.3.2)$$
Let $H\overline{\otimes}_A H\subset H\otimes_A H$ be the locus of coincidence of these
two structures, i.e., the space formed by elements $\sum_i u_{1i}\otimes u_{2i}$
satisfying
$$\sum_i u_{1i}a\otimes u_{2i} = \sum_i u_{1i}\otimes u_{2i}a\in H\otimes_A H, \quad \forall a\in A.
\leqno (A.3.3)$$

\vskip .2cm

\noindent {\bf (A.3.4) Example.} Let $A=k[X]$ be the coordinate algebra of
a smooth affine variety $X$ over $k$, and $H=\cD(X)$ be
the algebra of regular differential operators on $X$. Let $D: A\to A$ be a derivation,
so we consider  $D$ as an element of $H$. For $a\in A$ write $a'=D(a)$. Then
$D\otimes 1 + 1\otimes D\in H\overline{\otimes}_A H$. Indeed,
$$Da\otimes 1 + a\otimes D = aD\otimes 1 + a'\otimes 1 + a\otimes D,$$
while
$$D\otimes a + 1\otimes Da = D\otimes a + 1\otimes aD + 1\otimes a',$$
which is the same modulo  left $A$-linearity. 

\vskip .2cm

\proclaim (A.3.5) Lemma. (a) $H\btimes_A H$ is an associative algebra with respect to the
standard product
$$\biggl(\sum_i u_{1i}\otimes u_{2i}\biggr) \cdot \biggl( \sum_j v_{1j} \otimes v_{2j}\biggr) 
\quad =\quad \sum_{i,j} u_{1i} v_{1j}\otimes u_{2i}v_{2j}.$$
(b) $A$ is embedded as a subalgebra in $H\btimes_A H$, and the two $A$-module
structures on it are induced by the left and right multiplication with $A$.\hfill\break
(c) If $M,N$ are two left $H$-modules, then $M\otimes_A N$ is a left $H\btimes_A H$-module.
\qed

The following concept is a particular case of the concept of a bialgebroid as given in [Lu] and
[Xu]. We present here a self-contained exposition of this particular case. 
 A stronger structure,  involving the antipode and  called a Hopf algebroid, was defined by
J. Mr\v cun ([Mr], Def. 2.1).

\proclaim (A.3.6) Definition. Let $A$ be a commutative $k$-algebra. A left $A$-bialgebra
consists of:\hfill\break
(1) A possibly noncommutative algebra $H$ containing $A$. \hfill\break
(2) A morphism of left $A$-modules $\epsilon: H\to A$ which is a twisted ring homomoprhism:
$$\epsilon(uv) = \epsilon (u\cdot \epsilon(v)).$$
(3) A homomorphism of algebras $\Delta: H\to H\btimes_A H$ which is
identical on $A$, coassociative and satisfies the left and right counit properties with
respect to $\epsilon$. 

Note that the conditions in (2) imply that $A$ is a left $H$-module, via
$$(u,a)\mapsto u(a) = \epsilon (ua),\leqno (A.3.7)$$
and with respect to this structure $\epsilon$ is a homomoprhism of left $H$-modules. 
Further, if $\Delta(u)  = \sum_i u^{(1)}_i\otimes u^{(2)}_i$, then the counit property
implies that the right $A$-module structure on $H$ is expressible through
the left structure and $\Delta$:
$$ua = \sum_i u^{(1)}_i(a) \cdot u^{(2)}_i.\leqno (A.3.8)$$
Compare with Formula (6) of [Lu] and Formula (17) of [Xu]. Similarly, for the $H$-module 
structure on $A$, we have
$$u(ab) = \sum_i u^{(1)}_i(a) \cdot u^{(2)}_i(b).\leqno (A.3.9)$$
The last equality establishes part (b) of the following fact.

\proclaim (A.3.10) Proposition. Let $H$ be a left $A$-bialgebra, and $M,N$ be two left
$H$-modules. Then:\hfill\break
(a) $M\otimes_A N$ has a natural structure of a left $H$-module.\hfill\break
(b) If $M=A$, then the $H$-module structure on $A\otimes_A N=N$ is identical to the original
one on $N$. Similarly, if $N=A$. 

\proclaim (A.3.11) Theorem. Let $A$ be a commutative $k$-algebra, and $L\buildrel a\over\rightarrow
\Der (A)$ be a Lie-Rinehart $A$-algebra. Then:\hfill\break
(a)  The enveloping algebra $U_A(L)$ has a unique structure
of a left $A$-bialgebra, with $\epsilon$ given by (1.2.6) and $\Delta$ making each $x\in L$
primitive. \hfill\break
(b) If $L$ is projective over $A$, then all primitive elements of $U_A(L)$ lie in $L$.

\noindent {\sl Proof:} (a) This is Theorem 3.7 of [Xu]. The construction of $\Delta$
(without axiomatization of its properties) was explicitly given in [HS], \S 3.6. 
Part (b) is proved using the PBW theorem for $U_A(L)$, see [Ri], \S 3, similarly to
the well known fact for ordinary enveloping algebras ([Bo], \S II.5). \qed

\vskip .3cm

\noindent {\bf (A.4) Groupoid-like elements.} Let $A$ be a commutative $k$-algebra, and
$H$ a left $A$-bialgebra. For any commutative algebra $\Lambda$ and any homomorphism
$\xi: A\to\Lambda$ we have the change of scalars
$${}^\Lambda_\xi H = \Lambda\otimes_A H \leqno (A.4.1)$$
(using the left $A$-module structure on $H$ and the $A$-module structure on $\Lambda$
given by $\xi$). This is a $(\Lambda, A)$-bimodule. The counit $\epsilon$ of $H$
induces a moprhism of $\Lambda$-modules
$$  {}^\Lambda_\xi \epsilon: \,\,\,{}^\Lambda_\xi H \to\Lambda. 
\leqno (A.4.2)$$
Further, ${}^\Lambda_\xi H \otimes _\Lambda {}^\Lambda_\xi H$ has has one left $\Lambda$-module
structure and two right $A$-module structures, each commuting with the left one.
As in (A.3.3), we denote by ${}^\Lambda_\xi H \btimes_\Lambda {}^\Lambda_\xi H $
the locus of coincidence of these two $A$-module structures. Then the 
comultiplication $\Delta$ in $H$ induces a morphism of $(\Lambda, A)$-bimodules
$$\Delta: {}^\Lambda_\xi H \to \,\,\,{}^\Lambda_\xi H \btimes_\Lambda \, {}^\Lambda_\xi H,
\leqno (A.4.3)$$
which makes ${}^\Lambda_\xi H $ into a coassociative $\Lambda$-coalgebra with counit
${}^\Lambda_\xi  \epsilon$. Further, we get from (A.3.9) the following
relation between the left $\Lambda$- and the right $A$-module structures on
${}^\Lambda_\xi H$:
$$ua = \sum_i \,\, {}^\Lambda_\xi\epsilon (u^{(1)}_i a)\cdot u^{(2)}_i, \quad u\in
{}^\Lambda_\xi H, a\in A. \leqno (A.4.4)$$

\proclaim (A.4.5) Definition. An element $g\in \,\, {}^\Lambda_\xi H$ will be called
groupoid-like, if $g$ does not vanish anywhere on $\Spec(\Lambda)$, and
$\Delta(g)=g\otimes g$.

Compare [Mr], p. 270 for a related concept of a weakly group-like element. 

\vskip .2cm

Note that the counit property implies $\,\, {}^\Lambda_\xi \epsilon (g) \cdot g=g$, and so
${}^\Lambda_\xi \epsilon(g) = 1$. We will refer to $\xi$ as the target homomorphism of $g$,
and write $\xi=t(g)$. 

Next, we associate to $g$ another homomorphism $\eta: A\to\Lambda$, called the
source of $g$ and denoted $\eta = s(g)$, by setting
$$\eta(a) = \,\, {}^\Lambda_\xi \epsilon (ga).\leqno (A.4.6)$$
The groupoid-like property of $g$ implies that $\eta$ is indeed a ring
homomoprhism:
$$\eta(a) \eta(b) =  (\, {}^\Lambda_\xi \epsilon \otimes\,\,{}^\Lambda_\xi \epsilon)(ga\otimes gb) = 
(\, {}^\Lambda_\xi \epsilon \otimes \,\, {}^\Lambda_\xi \epsilon)((g\otimes g)\cdot ab) = $$
$$= (\,{}^\Lambda_\xi \epsilon \otimes \,\, {}^\Lambda_\xi \epsilon)(\Delta(g)\cdot ab) = \,\, 
{}^\Lambda_\xi \epsilon(ab) = \eta(ab).$$
We now consider the multiplication in $H$ as a morphism $m: H\otimes_k H\to H$.
It is  $A$-linear with respect to the $A$-module structure given by left
multiplication of $a\in A$ on the first factor in the source and on the target.
Therefore it descends to
$${}^\Lambda_\xi m: \,\,\,  {}^\Lambda_\xi H \otimes_k H = \Lambda\otimes_A H\otimes_k H
\to \Lambda\otimes _A H = \,\,  {}^\Lambda_\xi H. \leqno (A.4.7)$$
Let $g\in \,\,  {}^\Lambda_\xi H$ be a groupoid-like element with the
source homomorphism $\eta$.

\proclaim (A.4.8) Lemma. The map
$$\mu_g = \,\, {}^\Lambda_\xi m(g\otimes -): H\to  {}^\Lambda_\xi H, \quad h\mapsto 
\,\, {}^\Lambda_\xi m(g\otimes h),$$
is left $A$-linear with respect to the $A$-module structure on 
$ {}^\Lambda_\xi H$ given by $\eta: A\to\Lambda$ and by the $\Lambda$-module
structure on $ {}^\Lambda_\xi H$. 

\noindent {\sl Proof:} For $a\in A$ we haven by associativity of $m$,
$$\mu_g(av) = \,\, {}^\Lambda_\xi m(g\otimes av) = \,\, {}^\Lambda_\xi m(ga\otimes v) 
\buildrel (A.4.4)\over = \,\, {}^\lambda_\xi m(\,\, {}^\lambda_\xi \epsilon(ga)\cdot g\otimes v) =$$
$$= \,\, {}^\Lambda_\xi \epsilon (ga) \,\, {}^\Lambda_\xi m (g\otimes v) = \eta(a) \mu_g(v),$$
where the equality labelled (A.4.4) uses $\Delta(g)=g\otimes g$ and the formula (A.4.4). \qed

\vskip .2cm

By the lemma, the map $\mu_g$ descends to
$${}^\Lambda_\eta\mu_g: \,\, {}^\Lambda_\eta H\to \,\, {}^\Lambda_\xi H. \leqno (A.4.9)$$
By arguments similar to the above, we now arrive at the following fact.

\proclaim (A.4.10) Theorem. (a) If $h\in \,\, {}^\Lambda_\eta H$ is another groupoid-like
element with $t(h)=\eta$ and $s(h)=\zeta$, then
$$g*h := \,\, {}^\Lambda_\eta \mu_g (h) \in \,\, {}^\Lambda_\xi H$$
is a groupoid-like element with source $\zeta$ and target $\xi$. \hfill\break
(b) The operation $*$ defines a category $\underline{H}(\Lambda)$ with objects 
being algebra homomoprhisms $A\to\Lambda$ and $\Hom(\eta,\xi)$ being the set of
groupoid-like elements with source $\eta$ and target $\xi$.

\vskip .3cm

\noindent {\bf (A.5) Dualizing left bialgebras.} Let $A$ and $H$ be as before.

\proclaim (A.5.1) Definition. By an admissible filtration on $H$ we mean an increasing
filtration $H=\bigcup_{d=0}^\infty H_{\leq d}$ with $H_{\leq 0}=A$, compatible with
both algebra and coalgebra structures, i.e.,
$$H_{\leq d}\cdot H_{\leq d'}\subset H_{\leq d+d'}, \quad \Delta(H_{\leq d}) \subset
\sum_{d'+d''=d} H_{d'}\otimes H_{d''},$$
and such that each $H_{\leq d}/H_{\leq d-1}$ is projective of finite rank as a left $A$-module.

Assume that $H$ is equipped with an admissible filtration and, further, is {\em cocommutative}.
Then each
$$H_{\leq d}^\vee = \Hom_{A-\Mod}(H_{\leq d}, A)\leqno (A.5.2)$$
is a commutative algebra with unit $\epsilon$. We have two homomorphisms of algebras
$$\sigma_d, \tau_d: A\to H_{\leq d}^\vee, \quad \sigma_d(a)(u) = \epsilon (au) = a\epsilon(u),\quad
\tau_d(a)(u) = \epsilon (ua). \leqno (A.5.3)$$
The projective limit
$$H^\vee = \lim_{\longleftarrow} {}_d \, H_{\leq d}^\vee = \Hom_{A-\Mod}(H, A)\leqno (A.5.4)$$
is then a commutative topological algebra and we have two homomorphisms
$\sigma, \tau: A\to H^\vee$. Further, the multiplication in $H$ gives rise to comultiplication
$$H^\vee \to H^\vee\widehat{\otimes}_A H^\vee = \lim_{\longleftarrow} {}_d \,\,
H_{\leq d}^\vee \otimes_A H^\vee_{\leq d},\leqno (A.5.5)$$
with the $A$-module structures on the factors given by $\tau$ and $\sigma$. In other words,
$H^\vee$ is a topological commutative Hopf algebroid over $A$ in the sense of [Ra]. 
Thus the ind-scheme
$$\GG_1 = \Spf(H^\vee) = ``{\lim_{\longrightarrow}}'' {}_d \,\, \Spec(H^\vee_{\leq d})
\leqno (A.5.5)$$
is equipped with two morphisms $s,t: \GG_1\to\Spec(A)$ and the composition as in (A.1.2). 
In other terms, $\GG = (\GG_0=\Spec(A), \GG_1)$ is a category object in ind-schemes.
In particular, for each commutative algebra $\Lambda$ we have a category $\GG(\Lambda)$
whose objects are homomorphisms $A\to\Lambda$.

\proclaim (A.5.7) Proposition. The category $\GG(\Lambda)$ is identified
 with $\underline{H}(\Lambda)$.

\noindent {\sl Proof:} Consider a morphism $\Spec(\Lambda)\to\GG_1$, i.e., an algebra
homomorphism $f: H^\vee_{\leq d}\to\Lambda$ for some $d$. Combining it with $\tau$, we get
a homomorphism $\xi: A\to\Lambda$. Vieving $\Lambda$ as an $A$-module via $\xi$, we see
that $f$ is, in particular, a morphism of $A$-modules, i.e., an element of
$$\Hom_{A-\Mod}(\Hom_{A-\Mod}(H_{\leq d}, A), \Lambda) \simeq \Lambda\otimes_A H_{\leq d} = 
\,\, {}^\Lambda_\xi H_{\leq d}.$$
Let $g\in \,\,\Lambda\otimes_A H_{\leq d}$ be the element corresponding to $f$.
Then the fact that $f$ is a ring (and not just an $A$-module) homomorphism,
is translated into the fact that $g$ is groupoid-like. We leave the rest to the reader. \qed

\vskip .1cm

\proclaim (A.5.8) Definition. Let $L$ be a Lie-Rinehart $A$-algebra. An admissible filtration
in $L$ is an uncreasing filtration $L=\bigcup_{d=1}^\infty L_{\leq d}$ by $A$-modules,
compatible with the brackets and such that each $L_{\leq d}/L_{\leq d-1}$ is a projective
$A$-module of finite rank.

For example, if $L$ is itself projective of finite rank over $A$, we can take $L_{\leq 1}=L$. 

Given an admissible filtration on $L$, we define a filtration on $U=U_A(L)$ by
$$U_{\leq d} = {\operatorname{Span}}_A\biggl\{ x_{i_1} ... x_{i_p}, \quad x_{i_\nu}\in L_{\leq d_\nu}, 
\quad \sum_\nu d_\nu \leq d\biggr\}, \quad d\geq 0.\leqno (A.5.9)$$
This filtration is admissible in the sense of (A.5.1).

\proclaim (A.5.10) Proposition. Let $L$ be a Lie-Rinehart $A$-algebra with an admissible
filtration. Then: \hfill\break
(a) The multiplication in $U_A(L)$ makes $\GG = (\GG_0=\Spec(A), \GG_1)$
into a formal groupoid. We denote this groupoid by $e^L$.\hfill\break
(b) The identification
$$U_A(L) = \Hom^{\operatorname{cont}}_A(U_A(L)^\vee, A)$$
(continuous morphisms of modules) restricts to an identification
$$L = \Der^{\operatorname{cont}}_A (U_A(L)^\vee, A)$$
(continuous derivations).

Part (b) is an algebraic analog of the equality (A.1.3) and means that $L$ can be seen as
the Lie algebroid of the formal groupoid $e^L$.  Part (b) expresses the fact that $L$
coincides with the subalgebra of primitive elements in $U_A(L)$, see (A.3.11)(b).

To prove (a), note that we already have that $\GG$ is a category object in ind-schemes.
Further, $\GG_0=\Spec(A)$ is a scheme. Each $U_{\leq d}^\vee$ is a nilpotent
extension of $A$, so $\GG_1$ is a formal ind-scheme, and $s,t,e$ induce identity
on reduced subschemes. This latter fact also
 implies that for any $\Lambda$ all morphisms in $\GG(\Lambda)$ are
invertible. So $\GG$ is indeed a formal groupoid. 

\vskip .3cm

\noindent {\bf (A.6) Globalization.} The following rather restrictive concept
 will be sufficient for
our purposes.

\proclaim (A.6.1) Definition. A formal ind-scheme $\cX$ over $k$ will be called smooth, if
it is locally isomorphic to an ind-scheme of the form $Y\times \GD_n$
for some smooth algebraic variety $Y$ over $k$ and some $n\geq 0$. 

Thus, for a smooth formal ind-scheme $\cX$ we have the tangent sheaf $T_X$ which
is a locally free sheaf of $\cO_\cX$-modules of finite rank. In particular, we can speak about
Lie algebroids on $\cX$, as well as about $U(\cG)$, the enveloping algebra of a Lie algebroid
$\cG$, see (1.2). 
We also have the concept of an admissible filtration on $\cG$ by sheafifying (A.5.7).

Let $\cG$ be a Lie algebroid on $X$ equipped with an admissible filtration $\{\cG_{\leq d}\}_{d\geq 1}$. We
then define an algebra filtration $\{ U_{\leq d}(\cG)\}$ on $U(\cG)$ as in (A.5.8).
The considerations of (A.5) extend easily to this situation by working with $U(\cG)$
as a sheaf of topological left bialgebras over $\cO_{\cX}$. Thus, we define
$$\GA_d(\cG) = \shHom_{\cO_\cX}(U_{\leq d}(\cG), \cO_\cX). \leqno (A.6.2)$$
(Hom with respect to the left module structure.)
Then  $\GA_d(\cG)$ is a sheaf of  commutative $\cO_\cX$-algebras with two algebra
embeddings $\sigma,\tau: \cO_\cX\to\GA_d(\cG)$. 
 Let us further set 
$$\GA(\cG) = \lim_{\longleftarrow} {} _d \,\, \GA_d(\cG), \quad \GG_1 = \Spf \,\,\GA(\cG) = 
``\lim_{\longrightarrow} {''} {}_d \,\, \Spec  \, \GA_d(\cG). \leqno (A.6.3)$$
We have then, by globalizing (A.5.10):

\proclaim (A.6.4) Theorem.  The multiplication in $U(\cG)$ makes $\GG=(\GG_0=\cX, \GG_1)$
into a formal groupoid denoted $e^{\cG}$, and we have
$$\cG = \underline{\Der}^{\operatorname{cont}}_{\cO_\cX}(\GA(\cG), \cO_\cX) = e^* T_{\GG_1/\GG_0}.$$

Further, let $S$ be a scheme, and $y: S\to\cX$ be a moprhism. We have then a quasicoherent sheaf
$y^*U(\cG)$ of $\cO_S$-modules. By sheafifying (A.4.3), we see that $y^*U(\cG)$ is a sheaf
of unital $\cO_S$-coalgebras, so we can speak about its groupoid-like sections. Given any such
section $g$, we have, as in (A.4.6), another morphism $x: S\to\cX$ called the source of
$g$ (while $y$ is called the target). 

\proclaim (A.6.5) Theorem. Morphisms $S\to \GG_1$ are in bijection with pairs $(y,g)$
where $y: S\to \cX$ is a morphism, and $g$ is a groupoid-like section of $y^* U(\cG)$. 
\qed

\vskip .2cm

\noindent {\bf (A.6.6) Example.} Let $\cG=T_\cX$ with filtration
$\cG_{\leq d} = \cG$, $d\geq 1$. Then $U(\cG) = \cD_\cX$, and
$U_{\leq d}(\cG) = \cD^{\leq d}_X$ is the sheaf of differential
operators of order $\leq d$. The dual bundle of commutative algebras
$\GA_d$ is then identified with $J^d(\cO_\cX)$, the bundle of $d$-jets of
functions. So $\GA=J^\infty(\cO_\cX)$ is the bundle of infinite jets, which
is the same as $\widehat{\cO}_{\cX\times\cX, \cX}$, the formal completion
of $\cO_{\cX\times\cX}$ along the diagonal. Thus $\GG_1 = \cX^{(\infty)}_{\cX\times \cX}$
is the formal neighborhood of the diagonal in $\cX\times \cX$, with $s,t$ being the
projections and $e$ being the diagonal embedding. We will call $\GG=e^{T_\cX}$ the
{\em crystalline groupoid} of $\cX$. 

\vskip .3cm

\noindent {\bf (A.7) 
 Transitive Lie algebroids on the formal disk.}
Let $$0\to\cG^\circ\to\cG\buildrel\alpha\over\longrightarrow T_{\GD_n}\to 0\leqno (A.7.1)$$
be a transitive Lie algebroid on $\GD_n$. We assume that $\cG$ is equipped with an admissible
filtration $\{\cG_{\leq d}\}$.
Then  $\cG^\circ$ is is a bundle of Lie algebras (possibly of infinite rank)
over $\GD_n$, i.e., a Lie  $k[[x_1, ..., x_n]]$-algebra free as a $k[[x_1, ..., x_n]]$-module. We denote by
$$\Gg = \cG^\circ/(x_1, ..., x_n)\cG^\circ\leqno (A.7.2)$$
the fiber of $\cG^\circ$ at the origin of $\GD_n$. This is a Lie $k$-algebra. 

Let $\Mod(\cG)$ be the category of $\cG$-modules, i.e., of finite rank vector bundles
$\cE$ on $\GD_n$ with $\cG$-action, as in (1.1.4). Let also $\Rep(\Gg)$ be the category
of finite-dimensional representations of the Lie algebra $\Gg$ over $k$.

\proclaim (A.7.3) Theorem. The restriction
functor
$$\Res: \Mod(\cG)\to\Rep(\Gg), \quad  \cE\mapsto \cE/(x_1, ..., x_n)\cE$$
is an equivalence of categories. 

\vskip .1cm

\noindent {\bf (A.7.4) Example.} Let $\cG = T_{\GD_n}$. Then the theorem 
reduces to the well known fact 
that the category of finite rank  vector bundles on $\GD_n$ with flat connections
is equivalent to the category of finite-dimensional $k$-vector spaces.

\vskip .2cm

\noindent {\sl Proof of Theorem A.7.3:} Let $G=e^\Gg$ be the formal group integrating $\Gg$,
and $\GG=e^{\cG}$ be the formal groupoid integrating $\cG$. Thus $\GG_0=\GD_n$, and
$\GG_1=\Spf(\GA)$, where $\GA = \GA(\cG)$ is the topological Hopf algebroid constructed in (A.6.3).
Because $\cG$ is transitive, we have surjective homomoprhism 
$U(\cG)\to U(T_{\GD_n}) = \cD_{\GD_n}$. Thus the projection to the crystalline
groupoid of $\GD_n$,
$$p: \GG_1\to (\GD_n)_{\GD_n\times\GD_n}^{(\infty)} = \GD_n\times\GD_n,
\leqno (A.7.5)$$
makes $\GG_1$ into a $G$-bitorsor over $\GD_n\times\GD_n$, i.e., it has two $G$-actions,
 each of which makes it into a torsor. Let
$$\overline{\GA} = \GA/(x_1,..., x_n)\GA\leqno (A.7.6)$$
be the fiber of $\GA$ at $0\in\GD_n$ (with respect to the right $\cO_{\GD_n}$-module structure).
Then, with respect to the left module structure, $\overline{\GA}$ is still an
$\cO_{\GD_n}$-module. Geometrically,
$$\Spf(\overline{\GA}) = p^{-1}(\{0\}\times\GD_n). \leqno (A.7.7)$$
Restricting the identification of Proposition A.6.4(b), we get a homomorphism
$$\Gg\to\Der^{\cont}_{\cO_{\GD_n}}(\overline{\GA}),\leqno (A.7.8)$$
which is an isomoprhism, in virtue of $p$ being a torsor (hence isomorphic to
the projection $G\times\GD_n\times\GD_n\to\GD_n\times\GD_n$). In particular, 
$\Gg$ acts on $\overline{\GA}$ by derivations.

\vskip .2cm

Let now $V\in\Rep(\Gg)$. Set
$$\Ind(V) = (\overline{\GA}\otimes_k V)^{\Gg} = \shHom_{\cG^\circ}(U(\cG), V),
\leqno (A.7.9)$$
where $V$ is conisdered as a $\cG^\circ$-module via the projection
$\cG^\circ\to\Gg$. Then $\cG$ acts on $\Ind(V)$ because $\cG$ acts on
$\overline{\GA}$, the latter action being induced from the left $\cG$-action
on $\GA$. The natural identification of vector spaces
$$V\to \Res (\Ind(V))$$
is obvious. Let $\cE\in\Mod(\cG)$. Let us construct a natural isomoprhism of $\cG$-modules
$$\lambda_\cE: \cE\to \Ind(\Res(\cE)) = \shHom_{\cG^\circ}(U(\cG), \overline{\cE}).
\leqno (A.7.10)$$
Indeed, to construct $\lambda_\cE$, we just use the $U(\cG)$-action on $\cE$:
$$\lambda_{\cE}(\epsilon)(u) = u(\epsilon)(0), \quad \epsilon\in\cE, u\in U(\cG),$$
where the value at 0 means the image in $\overline{\cE}$. So we have constructed 
$\lambda_\cE$ as a morphism of $\cG$-modules. Further, after restriction to the
fibers at $0\in\GD_n$, the morphism $\lambda_\cE$ is an isomoprhism of $k$-vector
spaces. Now, the Nakayama lemma for $\GD_n=\Spf\,\, k[[x_1, ..., x_n]]$
implies that $\lambda_\cE$ is an isomorphism of vector bundles over $\cD_n$. 
\qed

\vskip 2cm

\centerline {\bf REFERENCES}

\vskip 1cm

\noindent [Bi] J.-M.  Bismut, {\em The Atiyah-Singer theorems: a probabilistic approach. I. The index theorem.}
 J. Funct. Anal. {\bf 57} (1984),  56-99. 

\vskip .1cm

\noindent [Bo] N. Bourbaki, {\em Groupes at Alg\'ebres de Lie}, 
Paris, Hermann, 1971.

\vskip .1cm

\noindent [CLP] J. M. Casas, M. Ladra, T. Pirashvili, {\em
Triple cohomology of Lie-Rinehart algerbas and the canonical class of
associative algebras}, preprint math.KT/0307354.

\vskip .1cm

 \noindent [C] K.~T.~Chen, {\em Iterated path integrals},  Bull. Amer. Math. Soc.  {\bf 83}  (1977), no. 5, 831--879. 

\vskip .1cm

\noindent [Co] A. Connes, {\em Noncommutative Geometry}, Academic Press, 1994.

 \vskip .1cm

\noindent [CDV] A. Connes, M. Dubois-Violette, {\em Yang-Mills algebra}, 
Lett. Math. Phys. {\bf 61}  (2002), 149-158, preprint mat.QA/0206205.

\vskip .1cm

\noindent [De] P. Deligne, {\em Le groupe fondamental de la droite projective
moins trois points}, in:  Galois groups over $Q$ (Berkeley, CA, 1987), 79-297, 
Math. Sci. Res. Inst. Publ., 16, Springer, New York, 1989.

\vskip .1cm

\noindent [Di] J. Dixmier, {\em Enveloping Algebras,} 
American Mathematical Society, Providence, RI, 1996. 

\vskip .1cm

\noindent [Ep] D. B. A. Epstein, {\em Natural tensors on Riemannian manifolds},
J. Differential Geom. {\bf 10} (1975),  631-645.

\vskip .1cm

\noindent [FS] B. Feigin, B. Shoikhet, {\em
On $[A,A]/[A,[A,A]]$ and on a $W_n$-action on the consecutive commutators of free associative algebra,}
preprint math.QA/0610410. 

\vskip .1cm

\noindent [GKZ] I.M. Gelfand, M.M. Kapranov, A.V. Zelevinsky,
{\em Discriminants, Resultants and Multidimensional Determinants},
Birkhauser, Boston, 1994. 

\vskip .1cm

\noindent [Gr] A. Grothendieck, J. Dieudonn\'e {\em EGA I}, Publ. Math. IHES, {\bf 4} (1960), 5-228.

\vskip .1cm

\noindent [Hab] W. Haboush, {\em Infinite-dimensional algebraic geometry: algebraic structures
 on $p$-adic groups and their homogeneous spaces,}
 Tohoku Math. J. (2){\bf  57} (2005), no. 1, 65-117.

\vskip .1cm

\noindent [HKK] K. A. Hardie, K. H. Kamps, R. W. Kieboom,
{\em A Homotopy 2-Groupoid of a Hausdorff Space}, Applied Cat. Structures,
{\bf 8} (2000), 209-234.

\vskip .1cm

\noindent [HS] V.~Hinich, V.~Schechtman, {\em Deformation theory and Lie algebra homology, I,II},
 Algebra Colloq. {\bf 4} (1997), no. 2, 213--240, 291--316.

\vskip .1cm

\noindent [Kal] R. Kallstrom, {\em Smooth modules over Lie algebroids I}, preprint math.AG/9808108.

\vskip .1cm

\noindent [Kap1] M. Kapranov, {\em Rozansky-Witten invariants via Atiyah classes},
Compositio Math. {\bf 115} (1999), 71-113. 

\vskip .1cm

\noindent [Kap2] M. Kapranov, {\em Noncommutative geometry and path
integrals}, preprint math.QA/0612411.

\vskip .1cm

\noindent [KM] M. Kapranov, Y.I. Manin, {\em Morita theory for operads},
Amer. J. Math. {\bf 123} (2001), 811-838. 

\vskip .1cm

\noindent [KV1] M. Kapranov, E. Vasserot, {\em Vertex algebras and the formal loop 
space}, Publ. Math. IHES, {\bf 100} (2004), 209-269. 

\vskip .1cm

\noindent [KV2] M. Kapranov, E. Vasserot, {\em Formal loops II: a  local Riemann-Roch theorem for
determinantal gerbes}, preprint math.AG/0509646, to appear in Ann. Sci. ENS.

\vskip .1cm

\noindent [Kas] M. Kashiwara, {\em Algebraic study of systems of partial differential equations
(Master's Thesis, Tokyo Univ. Dec. 1970)},
Memoires Soc. Math. France, {\bf 63} (1995), 1-72

\vskip .1cm

\noindent [Ko] S.  Kobayashi, {\em  La connexion des vari\'et\'es
 fibr\'ees II}, 
  C. R. Acad. Sci. Paris {\bf 238} (1954), 443-444.
\vskip .1cm

\noindent [Lu] J.-H. Lu, {\em Hopf algebroids and quantum groupoids},
Int. J. Math. {\bf 7} (1996), 47-70.

\vskip .1cm

 \noindent [Mack] K.~C.~H.~Mackenzie, {\em General Theory of Lie Groupoids and Lie Algebroids},
 London Mathematical Society Lecture Note Series, {\bf 213}. Cambridge University Press, Cambridge, 2005.

\vskip .1cm

\noindent [Macd] I. G. Macdonald, {\em Symmetric Functions and Hall Polynomials}, 
Oxford Univ. Press, 1995.

\vskip .1cm

\noindent [Man] Y. I. Manin, {\em Frobenius Manifolds, Quantum Cohomology and
Moduli Spaces}, Amer. Math. Soc. Providence RI, 1999.

\vskip .1cm

\noindent [Mr] J. Mr\v cun, {\em On duality between \'etale groupoids and Hopf algebroids},
J. Pure Appl. Alg. {\bf 210} (2007), 267-282. 

\vskip .1cm

\noindent [Ne] N. A. Nekrasov, {\em  Lectures on open strings, and noncommutative gauge fields},
preprint  hep-th/0203109. 

\vskip .1cm

\noindent [Po] A. M. Polyakov, {\em Gauge fields and space-time},
 Int. J. Mod. Phys. {A17S1} (2002), 119-136, preprint hep-th/0110196.

\vskip .1cm

\noindent [Ra] D. C. Ravenel, {\em Complex Cobordism and Stable Homotopy
Groups of Spheres}, Academic Press, 1986.  
 
\vskip .1cm

\noindent [Re] C. Reutenauer, {\em Free Lie Algebras,} Oxford Univ. Press, 
1993.

\vskip .1cm

\noindent [Ri] G. S. Rinehart, {\em  Differential forms on general commutative algebras,}
 Trans. Amer. Math. Soc. {\bf 108} (1963),  195-222.

\vskip .1cm

\noindent [Si] I.M. Singer, {\em On the master field in two dimensions}, 
in: 
Functional analysis on the eve of the 21st century, Vol. 1
 (New Brunswick, NJ, 1993), 263-281, 
Progr. Math., {\bf 131}, Birkhauser, Boston, MA, 1995. 

\vskip .1cm

\noindent [St] P. Stredder, {\em Natural differential operators on 
Riemannian manifolds and representations of the orthogonal and special
orthogonal group}, J. Diff. Geom. {\bf 10} (1975), 647-660. 

\vskip .1cm

\noindent [Xu] P. Xu, {\em Quantum groupoids}, Comm. Math. Phys. 
{\bf 216} (2001), 539-581. 

\vskip .1cm

\noindent [We] H. Weyl, {\em  Classical Groups,} Princeton Univ. Press, 1939.

\vskip 2cm

\end{document}